\documentclass{article}

\usepackage[square,comma, sort&compress]{natbib}
\usepackage[margin=1in]{geometry}
\usepackage{jmlr2e}
\usepackage{pgf}

\usepackage[utf8]{inputenc} 
\usepackage[T1]{fontenc}    
\usepackage{url}            
\usepackage{booktabs}       
\usepackage{amsfonts}       
\usepackage[protrusion=true,expansion,tracking=true]{microtype}  
\usepackage{nicefrac}
\usepackage{bbm}
\usepackage{xspace}

\usepackage{algpseudocode}
\usepackage{algorithm}

\usepackage{bm}
\usepackage[colorinlistoftodos]{todonotes}
\usepackage{amsmath}
\usepackage{amsthm}
\usepackage{amssymb}
\usepackage{thmtools}
\usepackage{thm-restate}
\usepackage{paralist}
\usepackage{comment}
\usepackage{float}
\usepackage{balance}
\usepackage{stmaryrd}
\usepackage{subcaption}
\usepackage{mathtools}

\usepackage{enumitem}
\usepackage{fancyhdr}
\usepackage{wrapfig}

\usepackage{multirow}
\usepackage{pifont}

\definecolor{darkblue}{RGB}{0,0,139}

\hypersetup{
  linkcolor=darkblue,    
  citecolor=darkblue,   
  urlcolor=darkblue       
}

\usepackage{tikz}
\usetikzlibrary{arrows.meta, positioning, shapes.geometric}

\newcommand{\innp}[1]{\left\langle #1 \right\rangle}

\newcommand*{\vsepfbox}[1]{%
  \begingroup
    \sbox0{\fbox{#1}}%
    \setlength{\fboxrule}{0pt}%
    \mbox{\kern-\fboxsep\fbox{\unhbox0}\kern-\fboxsep}%
  \endgroup
}

\newcommand{\dcprob}[1]{\Phi(#1)}
\newcommand{\dcasub}[2]{\hat{\Phi}_{#1}(#2)}
\newcommand{\dcgap}[1]{\Delta_{#1}}
\newcommand{\frankwolfejl}{\texttt{FrankWolfe.jl}\xspace}
\newcommand{\boscia}{\texttt{Boscia.jl}\xspace}
\newcommand{\lmo}{LMO\xspace}

\theoremstyle{plain} \numberwithin{equation}{section}
\newtheorem{theorem}{Theorem}[section]
\numberwithin{theorem}{section}

\newtheorem{fact}[theorem]{Fact}
\theoremstyle{definition}

\newtheorem{remark}[theorem]{Remark}

\graphicspath{{imgs/}}
\makeatletter
\def\mathcolor#1#{\@mathcolor{#1}}
\def\@mathcolor#1#2#3{%
  \protect\leavevmode
  \begingroup
    \color#1{#2}#3%
  \endgroup
}
\makeatother

\usepackage{amsmath}
\usepackage{amssymb}

\usepackage{charter}
\usepackage[vvarbb]{newtxmath}
\usepackage[capitalize,nameinlink,noabbrev]{cleveref}
\usepackage{tabularx}
\usepackage{booktabs}
\usepackage{multirow}
\usepackage{pifont}

\makeatletter
\newenvironment{proof*}[1][\proofname]{\par
  \pushQED{\qed}%
  \normalfont \partopsep=\z@skip \topsep=\z@skip
  \trivlist
  \item[\hskip\labelsep
        \itshape
    #1\@addpunct{.}]\ignorespaces
}{%
  \popQED\endtrivlist\@endpefalse
}
\makeatother

\title{Scalable DC Optimization via Adaptive Frank-Wolfe Algorithms}
\author{\name Sebastian Pokutta \email \href{mailto:pokutta@zib.de}{pokutta@zib.de} \\
\addr Institute of Mathematics, 
Technische Universität Berlin and \\ Zuse Institute Berlin, Germany
}
\begin{document}

\maketitle

\begin{abstract}
We consider the problem of minimizing a difference of (smooth) convex functions over a compact convex feasible region $P$, i.e., $\min_{x \in P} f(x) - g(x)$, with smooth $f$ and Lipschitz continuous $g$. 
This computational study builds upon and complements the framework of \citet{maskanRevisitingFrankWolfeStructured2025} by integrating advanced Frank-Wolfe variants to reduce computational overhead. We empirically show that constrained DC problems can be efficiently solved using a combination of the Blended Pairwise Conditional Gradients (BPCG) algorithm \citep{TTP2021} with warm-starting and the adaptive error bound from \citet{maskanRevisitingFrankWolfeStructured2025}. The result is a highly efficient and scalable projection-free algorithm for constrained DC optimization.
\end{abstract}

\section{Introduction}
In this paper, we focus on solving constrained Difference-of-convex (DC) optimization of the form
\begin{equation}
    \min_{x \in P} \dcprob{x} \doteq \min_{x \in P} f(x) - g(x),
\end{equation}
where $P$ is a compact convex feasible region, $f$ is an $L_f$-smooth convex function, and $g$ is a $L_g$-Lipschitz continuous convex function, following the setup in \citet{maskanRevisitingFrankWolfeStructured2025}. This problem falls into the class of DC programming problems, which are typically nonconvex but still somewhat well-structured optimization problems. We are in particular interested in the case where the feasible region $P$ is complicated, so that traditional projection-based methods are computationally prohibitive within the DCA algorithm and projection-free methods, in particular the Frank-Wolfe algorithm, are a natural choice given that they only rely on a linear minimization oracle to access $P$, which is often much cheaper than projections \citep{CP2021}. We require access to the gradients of $f$, subgradients of $g$ via (sub-)gradient oracles and access to $P$ via a linear minimization oracle.

\subsection*{Related Work}
Both DC Programming as well as Frank-Wolfe methods are very active fields of research. Here we will only be able to briefly discuss works that are most relevant to our setting here. 

\paragraph{DC Programming}

The difference-of-convex optimization framework provides a powerful mathematical structure for tackling a broad class of structured nonconvex optimization problems. The fundamental algorithmic approach for this problem class is the Difference-of-Convex Algorithm (DCA), first developed by \citet{tao1986algorithms} and over the last thirty years DCA and the DC framework have established themselves as particularly relevant in practical context, such as, e.g., indefinite quadratic optimization \citep{thi1997solving,an1998branch}, convex-concave minimax problems \citep{yuille2003concave,shen2016disciplined}, machine learning kernel selection \citep{argyriou2006dc}, hierarchical bilevel optimization \citep{hoai2009dc}, contextual decision making under uncertainty \citep{bennouna2024addressing}, domain adaptation with distribution discrepancy \citep{awasthi2024best}, structured sampling via determinantal processes \citep{mariet2015fixed}, adversarial robustness in deep learning \citep{awasthi2024dc}, and more recently in generalized polyhedral DC problems \citet{huong2024generalized}. An efficient practical improvement to the DCA algorithm is \emph{boosted DCA}, where the new iterate of the DCA does not simply replace the old one but a line search is performed, which can improve computational performance in many cases \citep{aragon2018boosted,aragonartachoAcceleratingDCAlgorithm2018,aragonArtachoBoostedDCAlgorithm2022,ferreiraBoostedDCAlgorithm2022}. This improvement via an additional line search can also be readily combined with the approach presented here, however in initial experiments we did not find significant additional improvements over the speed-up offered by BPCG with warm-starting and adaptive error bounds (see \cref{sec:boosted-dc}). An in-depth analysis, which is beyond the scope of this paper, is left for follow-up work. Also the theoretical understanding of DCA has evolved significantly over time. Initially, \citet{tao1997convex} established the asymptotic convergence properties, which were subsequently refined by \citet{lanckriet2009convergence} through simpler analyses under differentiability conditions; see also \citep{zhang2020complexity,kong2023cost,davis2022gradient}. A major breakthrough came with the finite-time convergence guarantees showing $\mathcal{O}(1/T)$ rates, independently developed by \citet{yurtsever2022cccp} and \citet{abbaszadehpeivasti2023rate}. Quite recently, \citet{rotaru2025tightanalysisdifferenceofconvexalgorithm} provided a tighter DCA convergence analysis with applications to the convergence rates of proximal methods. We refer the interested reader to \citet{le2018dc} for an in-depth discussion of the theoretical foundations, algorithmic developments, and practical applications.

\paragraph{Frank-Wolfe algorithms}

The Frank-Wolfe algorithm \citep{fw56}, also known as the conditional gradient method \citep{polyak66cg}, is a projection-free optimization method for constrained problems that has recently regained significant interest. This renewed attention is largely due to its suitability for problems with complicated feasible regions, where projection-based methods are often impractical. Originally designed for smooth convex optimization over polyhedral sets, the algorithm was later generalized by \citet{polyak66cg} to arbitrary compact convex constraint sets. The main appeal of Frank-Wolfe lies in its computational simplicity: it replaces potentially expensive or intractable projection steps with linear optimization over the constraint set. The modern revival of Frank-Wolfe methods can be traced to the influential work of \citet{jaggi13fw}, which sparked a wave of follow-up research. Since then, the Frank-Wolfe algorithm has been extended in numerous directions, including nonconvex settings \citep{lacoste2016convergence}, stochastic and block-coordinate variants \citep{osokin2016minding}, and various acceleration strategies such as away steps \citep{guelat1986some}, pairwise steps \citep{lacoste15}, in-face directions \citep{freund2017extended}, lazification \citep{BPZ2017,BPZ2017jour}, blending \citep{pok18bcg,CP2019,TTP2021}, sliding techniques \citep{lan2016conditional,LPZZ2017,qu2018non}, and Nesterov-style acceleration \citep{CDP2019,CDLP2021}. For a comprehensive overview of these developments, we refer the interested reader to \citet{CGFWSurvey2022}.

\paragraph{Frank-Wolfe for DC problems}

The connection between DC programming and Frank-Wolfe algorithms go back to at least \citet{Khamaru2019Convergence}. Building on this, \citet{yurtsever2022cccp} provided a strong relationship between DCA and Frank-Wolfe optimization via epigraphic reformulations. More recently, research has focused on developing Frank-Wolfe methods tailored to difference-of-convex optimization. For example, \citet{millan2023frank} proposed a direct Frank-Wolfe approach for DC problems, showing that all accumulation points of their algorithm are stationary points of the original problem. Furthermore, under additional weak-star convexity assumptions, they established $\mathcal{O}(1/T)$ convergence rates for both duality gap measures and objective function residuals. Finally, this computational study is based on \citet{maskanRevisitingFrankWolfeStructured2025} as we will see further below.

\subsection*{Contribution}

In this paper, we present a computational study that builds upon the framework of \citet{maskanRevisitingFrankWolfeStructured2025}. By integrating their adaptive error bound with advanced Frank-Wolfe variants, most notably the Blended Pairwise Conditional Gradients (BPCG) algorithm \citep{TTP2021} with warm-starting, we obtain a highly efficient and scalable projection-free approach for difference-of-convex (DC) optimization problems, allowing for the efficient solution of large-scale instances that were previously out of reach.

We demonstrate the superiority of this approach by extensive experiments on several problem classes, such as, e.g., difference-of-convex quadratics, complicated nonconvex functions, and Quadratic Assignment Problems (QAP). In all experiments, we observe substantial gains in computational efficiency, consistently outperforming baseline methods and often achieving significant reductions in \lmo calls (often by orders of magnitude) and runtime speedups exceeding factors of $100$. 

\section{Preliminaries and Notation}

In the following we will briefly recall some basic notions and facts that we will use later in the paper and we sometimes suggestively write $x_t \in P$ in our formulas to hint at the later use as iterates. 

\subsection{The Difference-of-Convex Algorithm (DCA)}
\label{sec:dca}

We want to solve the problem
\begin{equation}
    \label{eq:dc-problem}
    \min_{x \in P} \dcprob{x} = \min_{x \in P} f(x) - g(x),
\end{equation}
where $P$ is a compact convex set and $\dcprob{x} \doteq f(x) - g(x)$ is a difference of two convex functions with $f$ being $L$-smooth and $g$ being potentially non-smooth but admitting subgradients. For simplicity of exposition we will assume that $g$ is smooth as well and write $\nabla g$ for the gradient of $g$, however our arguments apply to the general case. In particular, the notion of gap below that uses $\nabla g$ applies also to any subgradient of $g$ in the non-smooth case. 

While $\dcprob{x}$ is a non-convex function, we can linearize the $g$-part of $\dcprob{x}$ at $x_t \in P$ to obtain a convex relaxation $\dcasub{t}{x}$ which is defined as
\begin{equation}
    \label{eq:dc-subproblem}
    \dcasub{t}{x} \doteq f(x) - g(x_t) - \innp{\nabla g(x_t), x - x_t},
\end{equation}
and which is an $L$-smooth and convex function and by the convexity of $g$ it holds
$$
\dcprob{x} \leq \dcasub{t}{x}. 
$$
From this inequality we obtain the notion of the \emph{difference-of-convex gap (DC gap)}
\begin{equation}
    \label{eq:dc-gap}
    \dcgap{x_t} \doteq \max_{x \in P} \dcprob{x_t} - \dcasub{t}{x} =  \max_{x \in P} f(x_t) - f(x) - \innp{\nabla g(x_t), x_t - x},
\end{equation}
and it holds that 
$$
\dcgap{x_t} = \max_{x \in P} \dcprob{x_t} - \dcasub{t}{x} \leq \max_{x \in P} \dcprob{x_t} - \dcprob{x},
$$
i.e., the DC gap is a \emph{lower bound on the primal gap}; for brevity we also write $\dcgap{t}$. Moreover, if we choose $x_{t+1} = \arg\max_{x \in P} \dcprob{x_t} - \dcasub{t}{x}$, then it holds that
\begin{equation}
    \label{eq:dc-progress}
    \dcgap{t} = \dcprob{x_t} - \dcasub{t}{x_{t+1}} \leq \dcprob{x_t} - \dcprob{x_{t+1}},
\end{equation}
i.e., the DC gap is also a \emph{lower bound on the primal progress} that can be achieved. As such, as long as the DC gap is not zero, we can make progress and together with standard arguments it follows that the DC gap is $0$ at $x$ if and only if $x$ is a stationary point; see e.g., \citet{millan2023frank}. Moreover, from \eqref{eq:dc-progress} we immediately obtain convergence of the DC gap via standard arguments: Add up \eqref{eq:dc-progress} for $t = 0, \dots, T$, telescope, and then divide by $T+1$ to obtain:
\begin{equation}
    \label{eq:dc-gap-convergence-simple}
    \frac{1}{T+1} \sum_{t=0}^T \dcgap{t} \leq \frac{1}{T+1} \sum_{t=0}^T \left(\dcprob{x_t} - \dcprob{x_{t+1}}\right) = \frac{\dcprob{x_0} - \dcprob{x_{T+1}}}{T+1} \leq \frac{\dcprob{x_0} - \dcprob{x^*}}{T+1},
\end{equation}
where the last inequality is by using that $x^*$ is a solution to \eqref{eq:dc-problem}, so that the average of the DC gaps (and hence its minimum!) converges to $0$ at a rate of $O(1/T)$. This argument is standard and can also be made inexact with additive error $\varepsilon$; see e.g., \citet{maskanRevisitingFrankWolfeStructured2025}, see also \cref{sec:adaptive-error-bound} for a similar reasoning. The argument also immediately gives rise to the DCA algorithm \citep{tao1986algorithms}, which is a popular method for solving \eqref{eq:dc-problem} and precisely uses the recursion as defined above, i.e.,
\begin{equation}
    \label{eq:dca-algorithm}
    \tag{DCA}
    x_{t+1} \gets \arg\min_{x \in P} \dcasub{t}{x}.
\end{equation}

We conclude this section with collecting useful properties from the arguments above.

\begin{fact}[Useful properties]
    \label{fact:useful-properties}
    Let $x_t \in P$ with $0 \leq t \leq T$ be a sequence of points (e.g., the iterates of the DCA algorithm) and $x^* = \arg\min_{x \in P} \dcprob{x}$ be a solution to \eqref{eq:dc-problem}. Then, it holds that
\begin{enumerate}
    \item $\dcprob{x} \leq \dcasub{t}{x}$ for all $x \in P$ and $\dcprob{x_t} = \dcasub{t}{x_t}$ by convexity of $g$, i.e., linearizing $g$ at $x_t$ provides an upper bound function for $\dcprob{x}$, which is exact at $x_t$ itself.
    \item $\dcgap{t} = \max_{x \in P} \dcprob{x_t} - \dcasub{t}{x} =  \max_{x \in P} f(x_t) - f(x) - \innp{\nabla g(x_t), x_t - x} = \dcprob{x_t} - \dcasub{t}{x_{t+1}}$, where $x_{t+1} \gets \arg\min_{x \in P} \dcasub{t}{x}$, i.e., the DCA algorithm picks $x_{t+1}$ to be the point that minimizes the linearization of $\dcprob{x}$ at $x_t$ and hence realizes the DC gap.
    \item $\dcgap{t} \leq \max_{x \in P} \dcprob{x_t} - \dcprob{x}$, i.e., the DC gap is a lower bound for the primal gap.
\end{enumerate}
\end{fact}

\subsection{The Frank-Wolfe Algorithm (FW)}

The Frank-Wolfe algorithm \citep{fw56}, also known as the conditional gradient method \citep{polyak66cg}, is a projection-free first-order method particularly well-suited for solving convex optimization problems over compact convex sets with potentially complicated constraints. In the context of DC programming, the Frank-Wolfe algorithm emerges as a natural choice for solving the convex subproblems \eqref{eq:dc-subproblem} that arise at each iteration of the DCA algorithm; see \cref{sec:dca-fw-revisited}. We will confine ourselves to the most relevant properties of the Frank-Wolfe algorithm and its vanilla variant here. However, more advanced variants exist as we will discuss later in \cref{sec:advanced-fw}.

\begin{algorithm}[h]
    \caption{Frank-Wolfe Algorithm \citep{fw56}}
    \label{alg:frank-wolfe}
    \begin{algorithmic}[1]
    \Require Initial point $x_0 \in P$, smooth convex objective function $h$
    \Ensure Iterates $x_1, x_2, \ldots \in P$
    \For{$t = 0, 1, 2, \ldots$}
    \State $v_t \gets \arg\min_{v \in P} \langle \nabla h(x_t), v \rangle$ \Comment{Linear minimization oracle}
    \State $\gamma_t \gets \frac{2}{t+2}$ \Comment{Step size (open loop)}
    \State $x_{t+1} \gets (1-\gamma_t) x_t + \gamma_t v_t$ \Comment{Convex combination}
    \EndFor
    \end{algorithmic}
    \end{algorithm}
    
The (vanilla) Frank-Wolfe algorithm for minimizing a smooth convex function $h(x)$ over a compact convex set $P$ is given in \cref{alg:frank-wolfe}. For a point $x \in P$, the \emph{Frank-Wolfe gap} is defined as
\begin{equation}
\label{eq:fw-gap}
G(x) \doteq \max_{v \in P} \langle \nabla h(x), x - v \rangle = \langle \nabla h(x), x - v^{FW} \rangle,
\end{equation}
where $v^{FW} \in \arg\min_{v \in P} \langle \nabla h(x), v \rangle$, which is often called a \emph{Frank-Wolfe vertex} and the function that provides such a solution for a given gradient is called a \emph{Linear Minimization Oracle (\lmo)}. The Frank-Wolfe gap has the important property that $G(x) = 0$ if and only if $x$ is an optimal solution to $\min_{x \in P} h(x)$. Moreover, it holds 
$$
h(x_t) - h(x^*) \leq \innp{\nabla h(x), x - x^*} \leq \max_{v \in P} \innp{\nabla h(x), x - v},
$$
where $x^* \in \arg\min_{x \in P} h(x)$, the first inequality is by convexity, and the second inequality is by maximality, making the Frank-Wolfe gap a natural stopping criterion and a measure of optimality.
The convergence analysis of the Frank-Wolfe algorithm relies on the smoothness of the objective function; see \citep{jaggi13fw,CGFWSurvey2022} for details. For an $L$-smooth convex function $h$ over a compact convex set $P$ with diameter $D = \max_{x,y \in P} \|x - y\|$, the Frank-Wolfe algorithm achieves a convergence rate of
\begin{equation}
\label{eq:fw-convergence}
h(x_t) - h(x^*) \leq \frac{2LD^2}{t+2} \qquad \text{ and } \qquad \min_{0 \leq t \leq T} G(x_t) \leq \frac{4 LD^2}{T+2}.
\end{equation}

\section{An adaptive error bound for DC subproblems}
\label{sec:adaptive-error-bound}

As shown in \citet{maskanRevisitingFrankWolfeStructured2025}, given a target accuracy $\varepsilon > 0$, we can allow for inexact solutions of the DCA subproblem, while ensuring convergence to a first-order stationary point of problem \eqref{eq:dc-problem}.

\begin{theorem}[Inexact DCA with convergence guarantee; \citet{maskanRevisitingFrankWolfeStructured2025}]
\label{thm:inexact-dca}
Suppose that the sequence $x_t$ is generated by an inexact-DCA algorithm designed to solve the subproblems described in \eqref{eq:dc-subproblem} approximately, in the sense that
$$
\dcasub{t}{x_{t+1}} - \dcasub{t}{x} \leq \frac{\varepsilon}{2}, \quad \forall x \in P,
$$
with $\varepsilon > 0$. Then, the following bound holds:
\begin{equation}
    \label{eq:dca-inexact-convergence}
    \min_{0 \leq t \leq T} \max_{x \in P} \left\{ f(x_t) - f(x) - \langle\nabla g(x_t), x_t - x \rangle \right\} \leq \frac{\dcprob{x_0} - \dcprob{x^*}}{T + 1} + \frac{\varepsilon}{2}.
\end{equation}
\end{theorem}

While \cref{thm:inexact-dca} allows for practical implementations of the DCA algorithm, it still has several drawbacks. One the one hand, we need to decide on $\varepsilon > 0$ ahead of time and on the other hand, we might solve the subproblems to a much higher target accuracy then actually needed for convergence, which in turn results in potentially excessive number of calls to the linear minimization oracle; our later numerical experiments confirm that this is actually happening. In a nutshell, \cref{thm:inexact-dca} establishes that solving the subproblems inexactly up to additive $\varepsilon > 0$, leads to a modification of the progress lower bound from \eqref{eq:dc-progress} to
\begin{equation}
    \label{eq:dc-progress-inexact}
    \underbrace{\dcprob{x_t} - \dcasub{t}{x_{t+1}}}_{\text{DC gap at $x_t$}} \leq \underbrace{\dcprob{x_t} - \dcprob{x_{t+1}}}_{\text{primal progress}} + \varepsilon,
\end{equation}
i.e., we have an additive error of $\varepsilon$ in the primal progress. In order to address these potential drawbacks \citet{maskanRevisitingFrankWolfeStructured2025} devised an adaptive error criterion (similar criteria have also been used before in \citet{TSPP2023,P2023}) that trades the additive error of $\varepsilon$ for a multiplicative error of a factor $1/2$ in the progress, via an arithmetic-geometric crossover. The resulting convergence rate is slower by a factor of $2$, however we can remove the additive error; compare \eqref{eq:dca-inexact-convergence} to \eqref{eq:dca-convergences}. In particular, this allows for solving to arbitrary precision simply by letting the algorithm ``run longer'', without the need for restarts or similar. We reproduce a full proof below for completeness. 

\begin{theorem}[Adaptive error bound and convergence rate for DCA; \citet{maskanRevisitingFrankWolfeStructured2025}]
    \label{thm:adaptive-dca}
Let $x_t \in P$ with $0 \leq t \leq T$ be the iterates of the DCA algorithm, where the subproblems are solved inexactly satisfying the condition:
\begin{equation}
    \label{eq:dca-error-bound}
    \dcasub{t}{x_{t+1}} \leq \dcprob{x_t} - \innp{\nabla \dcasub{t}{x_{t+1}}, x_{t+1} - v_{t+1}},
\end{equation}
with $v_{t+1} \in \arg\max_{v \in P} \innp{\nabla \dcasub{t}{x_{t+1}}, x_{t+1} - v}$, which is precisely a Frank-Wolfe vertex. Then, it holds that
\begin{equation}
    \label{eq:dca-convergences}
    \begin{aligned}
    \min_{0 \leq t \leq T} \dcgap{t} &= \min_{0 \leq t \leq T} \max_{x \in P} \dcprob{x_t} - \dcasub{t}{x} \\
    &= \min_{0 \leq t \leq T} \max_{x \in P} f(x_t) - f(x) - \innp{\nabla g(x_t), x_t - x} \\
    &\leq 2 \frac{\dcprob{x_0} - \dcprob{x^*}}{T+1}, 
    \end{aligned}
\end{equation}
where $x^* = \arg\min_{x \in P} \dcprob{x}$ and $x_{t+1}$ is the inexact solution to the subproblem, i.e., the iterates $x_t$ converge to a stationary point of the DC problem $\min_{x \in P} \dcprob{x}$ with a convergence rate of $O(1/T)$.
 
\begin{proof}
Assume that \eqref{eq:dca-error-bound} holds for all $t \in \{0, \ldots, T\}$. Then, we have that for all $x \in P$
\begin{equation}
    \label{eq:conv_transform}
\begin{aligned}
\dcasub{t}{x_{t+1}} &\leq \dcprob{x_t} - \innp{\nabla \dcasub{t}{x_{t+1}}, x_{t+1} - v_{t+1}} \\
& \leq \dcprob{x_t} - (\dcasub{t}{x_{t+1}} - \dcasub{t}{x}),
\end{aligned}
\end{equation}
where the second inequality is by convexity and maximality of $v_{t+1} \in P$, i.e., $\dcasub{t}{x_{t+1}} - \dcasub{t}{x} \leq \innp{\nabla \dcasub{t}{x_{t+1}}, x_{t+1} - x} \leq \innp{\nabla \dcasub{t}{x_{t+1}}, x_{t+1} - v_{t+1}}$ for all $x \in P$. We then rearrange \eqref{eq:conv_transform} to
$$
\dcasub{t}{x_{t+1}} \leq \frac{1}{2} \left(\dcprob{x_t} + \dcasub{t}{x}\right),
$$ 
and subtract $\dcprob{x_t}$ from both sides and rearrange to obtain for all $x \in P$
\begin{equation}
    \label{eq:progress-pusher}
    \begin{aligned}
    \frac{1}{2} \left(\dcprob{x_t} - \dcasub{t}{x}\right) &\leq \dcprob{x_t} - \dcasub{t}{x_{t+1}} \\
    &\leq \dcprob{x_t} - \dcprob{x_{t+1}},
    \end{aligned}
\end{equation}
where the last inequality is by convexity. In particular, for every $t$ we have 
\begin{equation}
    \label{eq:progress-pusher-max}
    \begin{aligned}
    \frac{1}{2} \dcgap{t} \doteq \frac{1}{2} \max_{x \in P} \left(\dcprob{x_t} - \dcasub{t}{x}\right) &\leq \dcprob{x_t} - \dcasub{t}{x_{t+1}} \\
    &\leq \dcprob{x_t} - \dcprob{x_{t+1}}.
    \end{aligned}
\end{equation}

We then sum up \eqref{eq:progress-pusher} over all $0 \leq t \leq T$, so that for all $x \in P$
\begin{equation}
    \label{eq:progress-sum}
    \sum_{t=0}^T \frac{1}{2} \dcgap{t} \leq \sum_{t=0}^T (\dcprob{x_t} - \dcprob{x_{t+1}}) = \dcprob{x_0} - \dcprob{x_{T+1}},
\end{equation}
and rearranging a final time results in 
\begin{equation}
    \label{eq:progress-sum-bound}
    \min_{0\leq t \leq T} \dcgap{t} \leq \frac{1}{T+1} \sum_{t=0}^T \dcgap{t} \leq 2 \frac{\dcprob{x_0} - \dcprob{x_{T+1}}}{T+1} \leq 2 \frac{\dcprob{x_0} - \dcprob{x^*}}{T+1}.
\end{equation}
\end{proof}        
\end{theorem}

Before we continue a few remarks are in order.

\begin{remark}[Efficiency of the subproblem error bound] 
    \label{rem:subproblem-error-bound}
First of all, note that we can replace the factor $1/2$ and $2$ by $1-\varepsilon$ and $1/(1-\varepsilon)$, respectively via a slight modification of the argument. Furthermore, while \cref{thm:adaptive-dca} ensures convergence without additive errors, an obvious question is, whether we are not solving the subproblem to much higher precision compared to \cref{thm:inexact-dca}. This is in fact not so, to see this we can rearrange \eqref{eq:dca-error-bound} to:
$$
\underbrace{\innp{\nabla \dcasub{t}{x_{t+1}}, x_{t+1} - v_t}}_{\text{Frank-Wolfe gap of $\dcasub{t}{x}$ at $x_{t+1}$}} \leq \underbrace{\dcprob{x_t} - \dcasub{t}{x_{t+1}}}_{\geq \dcgap{t}},
$$
i.e., we solve the subproblem $\min_{x \in P} \dcasub{t}{x}$ until its Frank-Wolfe gap drops below $\dcprob{x_t} - \dcasub{t}{x_{t+1}}$, which is at least as large as the DC gap at $x_t$; here $x_{t+1}$ denotes the iterate of the subsolver that it would return. In particular, for a fixed $\varepsilon > 0$, as long as $\dcgap{t} \geq \varepsilon$, we do not solve the subproblems more exactly than required in \cref{thm:inexact-dca}. We will see later though that in fact, \eqref{eq:dca-error-bound} is satisfied much earlier often requiring only a few iterations on the subproblem. Finally, the adaptive error bound criterion removes the need for deciding on $\varepsilon > 0$ ahead of time.
\end{remark}

\section{DCA and Frank-Wolfe}
\label{sec:dca-fw-revisited}

In the context of our DC subproblems \eqref{eq:dc-subproblem}, applying the Frank-Wolfe algorithm is a natural choice because the subproblems $\dcasub{t}{x}$ are $L$-smooth and convex by construction and $P$ is compact convex. Furthermore, as discussed in \cref{rem:subproblem-error-bound}, the Frank-Wolfe gap of the subproblem naturally appears in the adaptive error bound criterion \eqref{eq:dca-error-bound}, which allows us to easily link the required precision of the subproblem being solved with Frank-Wolfe variants to the outer DCA loop. This immediately lends itself to the algorithmic template presented in \cref{alg:adaptive-dca-fw}, which solves the DC subproblems $\min_{x \in P} \dcasub{t}{x}$ via some Frank-Wolfe algorithm $\mathcal A$ that uses as stopping criterion function a function $\tau_t$, i.e., 
\begin{equation}
\label{eq:fw-stopping-criterion}
\langle \nabla \dcasub{t}{y}, y - v \rangle \leq \tau_t(y) = \dcprob{x_t} - \dcasub{t}{y},
\end{equation}
where $v \in \arg\min_{v \in P} \langle \nabla \dcasub{t}{y}, v \rangle$ is the Frank-Wolfe vertex; we recover the variant of \cref{thm:inexact-dca}, simply by choosing $\tau_t = \varepsilon$. The algorithm can then be instantiated with different Frank-Wolfe algorithms such as vanilla Frank-Wolfe, Away-Step Frank-Wolfe, Pairwise Conditional Gradients, or Blended Pairwise Conditional Gradients for $\mathcal{A}$.

\begin{algorithm}[h]
\caption{Adaptive DCA with Frank-Wolfe}
\label{alg:adaptive-dca-fw}
\begin{algorithmic}[1]
\Require Initial point $x_0 \in P$, DC objective $\dcprob{x} = f(x) - g(x)$, Frank-Wolfe variant $\mathcal{A}$
\Ensure Sequence of iterates $x_1, x_2, \ldots \in P$
\For{$t = 0, 1, 2, \ldots$}
\State Define subproblem: $\dcasub{t}{x} = f(x) - g(x_t) - \langle \nabla g(x_t), x - x_t \rangle$
\State Define stopping criterion: $\tau_t(y) = \dcprob{x_t} - \dcasub{t}{y}$ for $y \in P$ \Comment{Adaptive error bound \eqref{eq:dca-error-bound}}
\State $x_{t+1} \gets \mathcal{A}(\dcasub{t}{\cdot}, x_t, \tau_t)$ \Comment{Call Frank-Wolfe variant}
\EndFor
\end{algorithmic}
\end{algorithm}

With $\tau_t = \varepsilon$ and $\mathcal A$ being the vanilla Frank-Wolfe algorithm (with some step-size strategy), we obtain the version in \citet{maskanRevisitingFrankWolfeStructured2025}, who unrolled the Frank-Wolfe algorithm and baked it directly into the DCA iterations. In contrast, we kept the modular structure to vary the subsolver easily.

\subsection{Leveraging advanced FW algorithms}
\label{sec:advanced-fw}

While \citet{maskanRevisitingFrankWolfeStructured2025} and \citet{millan2023frank} considered primarily the vanilla Frank-Wolfe algorithm \citep{fw56}, the former mentioning improved convergence rates in some special cases, many advanced variants of the base algorithm exist. These variants, including Away-Step FW (AFW) \citep{guelat1986some}, Fully-Corrective Frank-Wolfe (FCFW) \citep{holloway1974extension}, Decomposition-Invariant Conditional Gradients (DICG) \citep{garber2016linear,bashiri2017decomposition}, Pairwise Conditional Gradients (PCG) \citep{lacoste15}, Blended Conditional Gradients (BCG) \citep{BPTW2018}, Boosted Frank-Wolfe (BFW) \citep{CP2020boost}, Blended Pairwise Conditional Gradients (BPCG) \citep{TTP2021} often exhibit up to linear convergence rates depending on feasible region and function properties. Moreover, they provide additional benefits such as lazification \citep{LPZZ2017,BPZ2017,BPZ2017jour} and blending \citep{BPTW2018,CP2019,TTP2021}, where the former allows to aggressively reuse information from previous \lmo calls and the latter blending in different types of descent steps for improved convergence properties. Finally, active set-based algorithms (including AFW, PCG, BCG, BFW, and BPCG), i.e., those that explicitly maintain the current iterate as convex combinations allow for warm-starting, which in turn can significantly reduce computational overhead when called repeatedly as heavily exploited in conditional gradient-based MINLP solver \boscia \citep{HTBP2022}. Moreover, there exist also many modern, highly-efficient line search strategies specifically for Frank-Wolfe methods \citep{pedregosa2018step,P2023,HBMP2025} that go beyond the original open-loop step-size and the short-step strategy; we refer the reader to \citet{CGFWSurvey2022} for an in-depth discussion.

Based on the benchmarks in \citet{BDHKPTVW2024}, we opted for using the Blended Pairwise Conditional Gradients (BPCG) algorithm of \citet{TTP2021}. While it is slightly slower than the Blended Conditional Gradients algorithm of \citet{BPTW2018}, it offers superior numerical stability (heavily exploited in \boscia) due to blending with pairwise steps rather than projected gradient descent steps. BPCG allows for warm-starting by reusing the active set of a prior call and is adaptable to well-known function properties and properties of the feasible region. In particular, it achieves accelerated sublinear (for the sharp objectives) to linear (for strongly convex objectives) convergence rates over polytopes \citep{kerdreux2019restarting,KDP2018jour} (see also \citet{CGFWSurvey2022}). Similarly, BPCG achieves improved convergence rates over uniformly convex sets \citep{kerdreux2021projection} and strongly convex sets \citep{garber2015faster}; also for the local case \citep{KAP2021}. 

We also tested other variants including Away-Step Frank-Wolfe and Pairwise Conditional Gradients, Blended Conditional Gradients, as well as the vanilla Frank-Wolfe with tracking the active set for warm-starting as these are readily available in \frankwolfejl. However, their performance was significantly inferior to BPCG. Moreover, we note that BPCG offers further potential improvements to reduce computational overheads, such as e.g., lazification, which potentially reduces the number of \lmo calls even further, however BPCG with warm-starting and the adaptive error bound already reduces the overall number or required \lmo calls so dramatically that we opted for the tighter gap control offered by the non-lazified variant. Finally, we mention that BPCG can also be combined with active set corrections from \citet{HRBDP2025}, basically one-shot-ing the solution of the subproblem over the current active set via an auxiliary solver. This is applicable only to the case when $\min_{x \in P} \dcasub{t}{x}$ is a convex quadratic problem. 

\section{Numerical experiments}

We now present our numerical experiments, demonstrating the effectiveness of the combination of the adaptive error bound criterion with warm-starting. All experiments are conducted in Julia 1.11.5 and \cref{alg:adaptive-dca-fw} has been implemented in the \frankwolfejl package \citep{BCP2021,BDHKPTVW2024}, so that we can directly use the latest Frank-Wolfe algorithms and line search strategies. Our DCA code will be contributed to the \frankwolfejl package and the code to run the computational experiments will be made separately available on GitHub. All tests were performed on a Mac Pro M1 Max. For visualization purposes, we (additively) shift the primal objective values by their minimal value to ensure non-negative values for logarithmic plotting. The early stopping via the adaptive error bound \eqref{eq:dca-error-bound} is realized via the callback mechanism available in \frankwolfejl, removing the need for restarts or halving-strategies as done in \citet{maskanRevisitingFrankWolfeStructured2025}. All experiments are run with a single thread except for the larger Quadratic Assignment Problem (QAP) instances where \frankwolfejl employs the Hungarian method from the \texttt{Hungarian.jl} Julia package within the \lmo of the Birkhoff polytope which internally is multi-threaded; this only affects absolute timing for the large QAP instances but not relative differences. Note, that still some large QAP instances require several weeks of compute time for the slower algorithms for achieving the, admittedly relatively tight, target precision.

\subsection{Experimental setup}
\label{sec:experiments}

We evaluate DCA solver variants based on \cref{alg:adaptive-dca-fw} with different Frank-Wolfe variants and stopping criteria:

\begin{itemize}
    \item \textbf{DCA-FW}: Standard Frank-Wolfe with fixed $\varepsilon$ for subproblems as required by \cref{thm:inexact-dca}; basically the algorithm in \citet{maskanRevisitingFrankWolfeStructured2025} except for employing the much faster secant line search strategy from \citet{HBMP2025}. Note, that we also did preliminary tests with both the agnostic $\gamma_t = \frac{2}{t+2}$ step-size strategy as well as short-steps, which maximize progress over the smoothness approximation (see \citet{demyanov1970approximate,CGFWSurvey2022}), however both were strictly inferior to the secant line search strategy.
    \item \textbf{DCA-FW-ES}: Standard Frank-Wolfe with adaptive stopping criterion \eqref{eq:dca-error-bound} as required by \cref{thm:adaptive-dca} for early stopping.
    \item \textbf{DCA-BPCG}: Blended Pairwise Conditional Gradients of \citet{TTP2021} with fixed $\varepsilon$ for subproblems as required by \cref{thm:inexact-dca}.
    \item \textbf{DCA-BPCG-ES}: Blended Pairwise Conditional Gradients with adaptive stopping criterion \eqref{eq:dca-error-bound} as required by \cref{thm:adaptive-dca} for early stopping.
    \item \textbf{DCA-BPCG-WS}: Blended Pairwise Conditional Gradients with warm-starting and fixed $\varepsilon$ for subproblems as required by \cref{thm:inexact-dca}.
    \item \textbf{DCA-BPCG-WS-ES}: Blended Pairwise Conditional Gradients with warm-starting and adaptive stopping criterion \eqref{eq:dca-error-bound} as required by \cref{thm:adaptive-dca} for early stopping.
\end{itemize}

All algorithms use the secant line search strategy from \citet{HBMP2025}, which is currently the fastest line search strategy for Frank-Wolfe algorithms and the default strategy in \frankwolfejl. In particular for quadratics this line search strategy terminates with as an optimal solution with two function evaluations. For all experiments, we impose a maximum of number of DCA iterations (between $200$ and $500$ depending on the experiment) and we stop with a DCA gap less than $10^{-6}$, which may or may not be achieved within the DCA iteration limit. We also impose an iteration limit on the inner iterations of the Frank-Wolfe variant (between $10000$ and $50000$ depending on the experiment) and we stop with a Frank-Wolfe gap less than $10^{-6}/2$ or when the adaptive error bound is satisfied, for those variants using it. The feasible regions vary by experiment but include the $k$-sparse polytope, probability simplex, $\ell_1$-ball, and the Birkhoff polytope as specified below; see \citet{BDHKPTVW2024,BCP2021} for detailed descriptions of these constraint sets. To assess scalability, we test various instance sizes ranging from small instances with $n=10$ to large instances with roughly $n=1000$. For each instance size, we generate multiple random instances unless stated otherwise. We also observed that the algorithms do not necessarily converge to the same primal solution, which is not unexpected as the considered problems are non-convex and the individual trajectories of the algorithms can differ vastly; see \cref{sec:additional-numerical-experiments} for some individual runs.

We use performance profiles for various performance measures (DCA iterations, time, and \lmo calls) to summarize the performance of the different algorithmic variants. The performance profiles are log-log plots, and the convergence plots in \cref{sec:additional-numerical-experiments} use a logarithmic scale on the $y$-axis, and also on the $x$-axis when reporting \lmo calls. For the performance profiles, an instance is considered solved if the DCA gap is less than $10^{-6}$, unless stated otherwise. 

Note, that at first one might assume that all algorithms are equivalent in terms of DCA iterations. While this would be true if we were able to run all subsolvers to sufficient accuracy (and assuming that they are trajectory equivalent, which may or may not be the case), this is not realistic, as it would be prohibitively expensive for all but the fastest algorithms. Therefore, we have set iteration limits for the subsolvers in our tests to prevent stalling; however, this means that the subsolver solutions might be of lower accuracy. This is precisely where the better performance of the BPCG variants becomes apparent: they achieve significantly higher convergence speed, providing sufficiently accurate solutions within the iteration limit.

\subsection{Summary of the experiments}

Before considering the specific experiments, here we summarize our key findings. Note that these findings are very consistent across the experiments. 

\paragraph{Adaptive error bound efficiency.} The adaptive error bound criterion proves to be highly efficient when used for early stopping (ES). Effectively it dynamically adjusts the subproblem solution precision based on the current DC gap, avoiding unnecessary computational overhead in early iterations. This can be clearly seen in the instance plots' number of required \lmo calls for the subproblems in the early iterations; see \cref{sec:additional-numerical-experiments}.

\paragraph{Warm-starting benefits.} Warm-starting (WS) reuses both the last iterate $x_t$ and its associated active set, which contains points returned from previous linear minimization oracle (\lmo) calls which then can be reused. This strategy is particularly beneficial towards the end of the optimization process when the algorithm is close to the optimal solution, so that the subproblems change only slightly and warm-starting significantly accelerates convergence in this regime.

\paragraph{Iteration count trade-offs.} The ES solvers sometimes exhibit slightly slower convergence in terms of DCA iterations, which aligns with the theoretical prediction in \cref{thm:adaptive-dca}, where we trade a factor of 2 in the convergence rate for the removal of additive error terms.

\paragraph{\lmo call efficiency.} ES solvers require significantly fewer \lmo calls in the beginning phases of optimization. However, as higher accuracy is required in later iterations, the number of \lmo calls naturally increases to meet the adaptive stopping criterion. Similarly, the WS solvers require significantly fewer \lmo calls in the end, as warm-starting cuts out the initial phase when subproblems are very similar. 

\begin{figure}[htbp]
    \centering
        \includegraphics[width=0.99\textwidth]{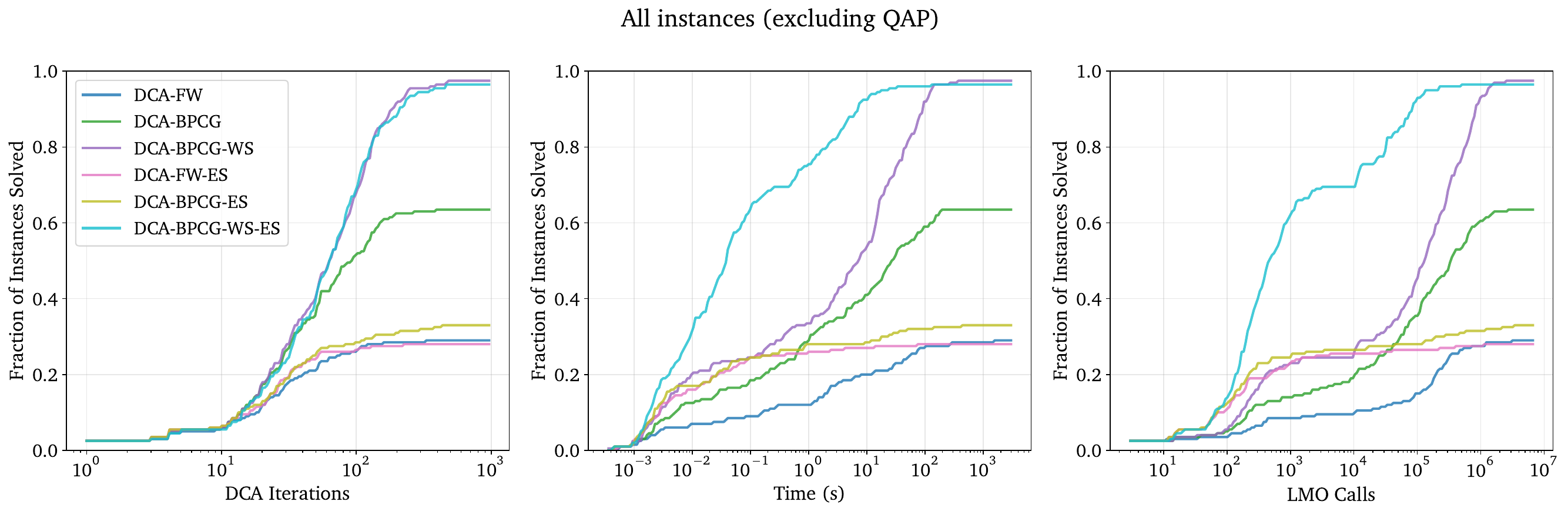}
    \caption{Performance profile over all tested instances (excluding QAP) as described in \cref{sec:experiments}. The plot demonstrates the superior efficiency of the DCA-BPCG-WS-ES variant (combining warm-starting and early stopping) across the entire range of problem instances.}
    \label{fig:performance-profile-all}
\end{figure}

\paragraph{Combination of adaptive error bound together with warm-starting.} While each of the two strategies ES or WS are very beneficial in terms of reducing computational cost, the combination of the two strategies is the most efficient by a huge margin, enabling multiple orders of magnitude fewer \lmo calls and reducing running times by a similar factor as shown in \cref{fig:performance-profile-all}. In particular, the combination of the two strategies enables us to tackle instance sizes that would otherwise be out of reach. In fact, in many cases, the improvement is so substantial that we cannot precisely quantify it because the baseline solvers fail to complete within the time limit.

Results across all considered instances (except for QAP) are summarized in \cref{tab:dca_fw_statistics} and \cref{fig:performance-profile-all}. The table uses shifted geometric means as customary; more results and examples of representative individual runs are to be found in \cref{sec:additional-numerical-experiments}. We excluded QAP from this overview as the instances typically do not solve to optimality and we thus we report these separately in detail in \cref{sec:qap}. As can be seen, the DCA-BPCG-WS-ES variant (combining warm-starting and early stopping) is the most efficient both in time and \lmo calls being often by multiple orders of magnitude faster (geometric means) than the other variants. For a very small number of experiments ($\approx 1\%$) it seems that DCA-BPCG-WS outperform DCA-BPCG-WS-ES as can be seen that the very top of the performance profile chart. However this is mostly an artifact of our experimental design and the hard cap on DCA iterations. By \cref{thm:adaptive-dca} we know that the primal progress may be slower in DCA iterations when using the adaptive error bound criterion while each of the DCA iterations is much faster. That is precisely what we see here: these select few experiments finish via DCA-BPCG-WS very close to the DCA iteration limit and the DCA-BPCG-WS-ES variant would require a few extra iterations to close the gap. This can be also seen when compared to the required compute time: DCA-BPCG-WS-ES requires several orders of magnitude less time (geometric mean) than DCA-BPCG-WS so that a few extra iterations are not a problem.

\begin{table}[htbp]
\centering
\caption{Performance Statistics over all instances (excluding QAP; see \cref{sec:qap}) [Shifted Geometric Means, shift: $1.0$].}
\label{tab:dca_fw_statistics}
\resizebox{\textwidth}{!}{
\begin{tabular}{lrrrrrrrrrrrrrrrrrr}
\toprule
Size & \multicolumn{3}{c}{DCA-FW} & \multicolumn{3}{c}{DCA-BPCG} & \multicolumn{3}{c}{DCA-BPCG-WS} & \multicolumn{3}{c}{DCA-FW-ES} & \multicolumn{3}{c}{DCA-BPCG-ES} & \multicolumn{3}{c}{DCA-BPCG-WS-ES} \\
 & Iter & Time & LMO & Iter & Time & LMO & Iter & Time & LMO & Iter & Time & LMO & Iter & Time & LMO & Iter & Time & LMO \\
\midrule
$n=10$ (10) & 24.65 & 1.57 & 6152.00 & \textbf{8.90} & 0.06 & 106.19 & \textbf{8.90} & 0.04 & 114.97 & 20.05 & 0.76 & 434.32 & 16.44 & 0.85 & 483.16 & 9.65 & \textbf{3.04e-03} & \textbf{42.05} \\
$n=20$ (10) & 63.91 & 5.20 & 4.98e+05 & \textbf{16.06} & 0.16 & 3361.08 & 16.24 & 0.16 & 1544.55 & 71.32 & 4.87 & 7.82e+04 & 71.15 & 5.20 & 6.41e+04 & 16.50 & \textbf{2.83e-03} & \textbf{90.21} \\
$n=30$ (10) & 61.92 & 4.32 & 1.06e+05 & 21.09 & 0.28 & 2861.19 & 21.26 & 0.37 & 2673.25 & 29.04 & 0.81 & 1056.74 & 28.95 & 1.14 & 705.45 & \textbf{19.64} & \textbf{6.38e-03} & \textbf{154.21} \\
$n=50$ (10) & 67.22 & 8.51 & 8.19e+04 & 34.84 & 2.15 & 1.36e+04 & \textbf{22.69} & 0.75 & 3501.54 & 50.46 & 5.15 & 7140.91 & 49.67 & 4.52 & 5760.02 & 30.31 & \textbf{0.07} & \textbf{330.21} \\
$n=100$ (20) & 144.57 & 29.09 & 8.85e+05 & 44.70 & 1.89 & 3.83e+04 & 41.49 & 1.15 & 6279.90 & 118.36 & 15.92 & 1.22e+05 & 125.41 & 17.24 & 1.24e+05 & \textbf{36.18} & \textbf{0.09} & \textbf{284.85} \\
$n=150$ (10) & 60.89 & 28.57 & 1.13e+05 & 39.72 & 7.74 & 3.90e+04 & \textbf{26.36} & 2.94 & 6116.85 & 54.70 & 13.09 & 1.05e+04 & 50.80 & 11.65 & 9254.55 & 27.58 & \textbf{0.02} & \textbf{135.35} \\
$n=200$ (20) & 155.52 & 78.73 & 1.39e+06 & 81.08 & 12.64 & 3.01e+05 & \textbf{58.71} & 3.38 & 5.00e+04 & 175.92 & 42.56 & 3.38e+05 & 165.38 & 34.04 & 2.33e+05 & 64.16 & \textbf{0.13} & \textbf{711.14} \\
$n=300$ (20) & 221.76 & 212.14 & 2.16e+06 & 123.29 & 37.87 & 6.14e+05 & \textbf{81.22} & 12.88 & 2.24e+05 & 246.98 & 123.15 & 9.43e+05 & 261.51 & 122.23 & 9.53e+05 & 90.48 & \textbf{0.57} & \textbf{2106.04} \\
$n=400$ (20) & 259.76 & 292.36 & 1.74e+06 & 218.06 & 93.76 & 1.05e+06 & 91.87 & 15.53 & 1.36e+05 & 279.67 & 228.21 & 1.54e+06 & 260.28 & 146.31 & 9.45e+05 & \textbf{80.69} & \textbf{0.63} & \textbf{1876.60} \\
$n=500$ (20) & 180.95 & 291.42 & 1.14e+06 & 131.41 & 79.50 & 5.66e+05 & \textbf{70.59} & 21.53 & 1.68e+05 & 223.16 & 180.22 & 7.62e+05 & 206.73 & 132.95 & 5.74e+05 & 87.64 & \textbf{2.34} & \textbf{8266.83} \\
$n=600$ (10) & 400.69 & 408.89 & 2.33e+06 & 217.09 & 135.01 & 1.70e+06 & \textbf{121.52} & 30.45 & 3.78e+05 & 399.20 & 306.11 & 1.86e+06 & 402.62 & 222.93 & 1.39e+06 & 135.39 & \textbf{1.35} & \textbf{7345.49} \\
$n=700$ (10) & 421.90 & 298.50 & 1.92e+06 & 298.63 & 109.60 & 9.26e+05 & 92.27 & 12.81 & 9.15e+04 & 399.95 & 260.05 & 1.74e+06 & 306.37 & 127.02 & 9.67e+05 & \textbf{80.22} & \textbf{1.14} & \textbf{3409.07} \\
$n=800$ (10) & 365.52 & 579.95 & 3.42e+06 & 356.54 & 346.69 & 3.09e+06 & 128.18 & 40.42 & 3.97e+05 & 368.72 & 251.60 & 1.38e+06 & 360.66 & 209.15 & 1.35e+06 & \textbf{125.44} & \textbf{1.66} & \textbf{5473.48} \\
$n=900$ (10) & 221.82 & 331.17 & 8.10e+05 & 159.99 & 191.99 & 6.21e+05 & 88.35 & 46.95 & 1.90e+05 & 200.90 & 121.79 & 2.68e+05 & 170.29 & 87.36 & 2.10e+05 & \textbf{86.70} & \textbf{2.21} & \textbf{5270.98} \\
$n=1000$ (10) & 236.52 & 239.21 & 4.59e+05 & 151.12 & 90.62 & 2.32e+05 & 88.63 & 27.89 & 7.70e+04 & 228.92 & 204.13 & 3.94e+05 & 197.63 & 129.56 & 2.79e+05 & \textbf{88.62} & \textbf{1.92} & \textbf{3084.11} \\
\bottomrule
\end{tabular}
}
\end{table}

\subsection{Difference of Convex Quadratic Functions}
\label{sec:diff-quadratics}

\begin{figure}[htbp]
    \centering
        \includegraphics[width=0.99\textwidth]{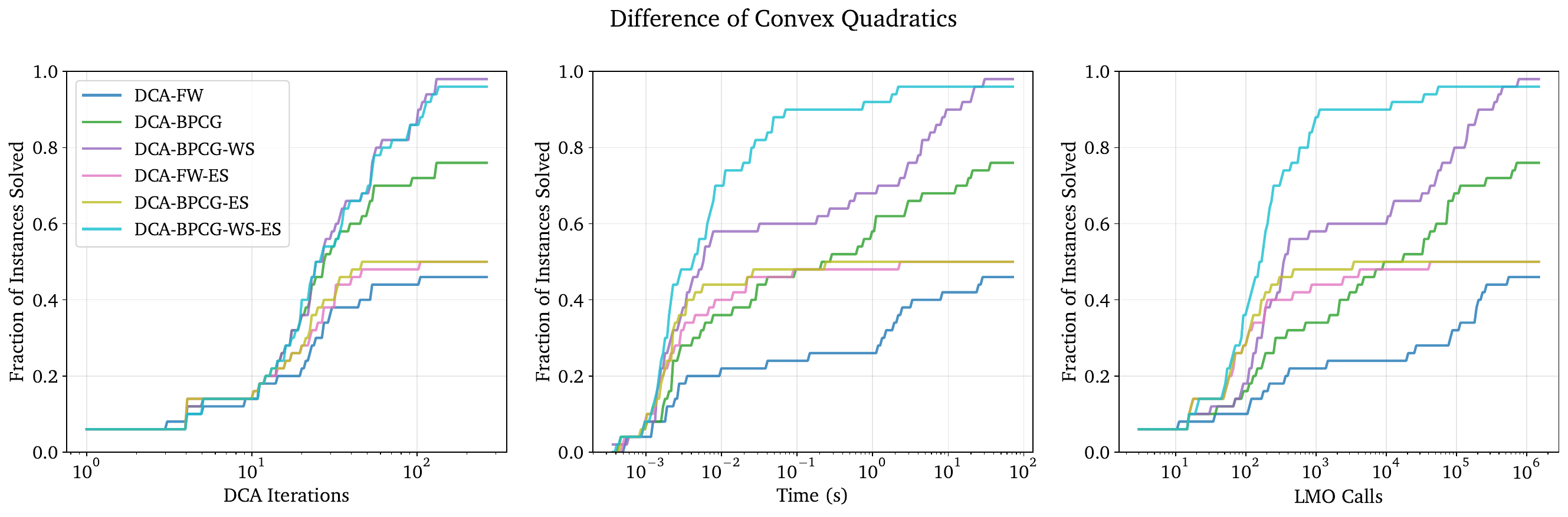}
    \caption{Performance profile over small to medium-sized instances of differences of convex quadratics from \cref{sec:diff-quadratics}.}
    \label{fig:performance-profile-simple}
\end{figure}

\begin{figure}[htbp]
    \centering
        \includegraphics[width=0.99\textwidth]{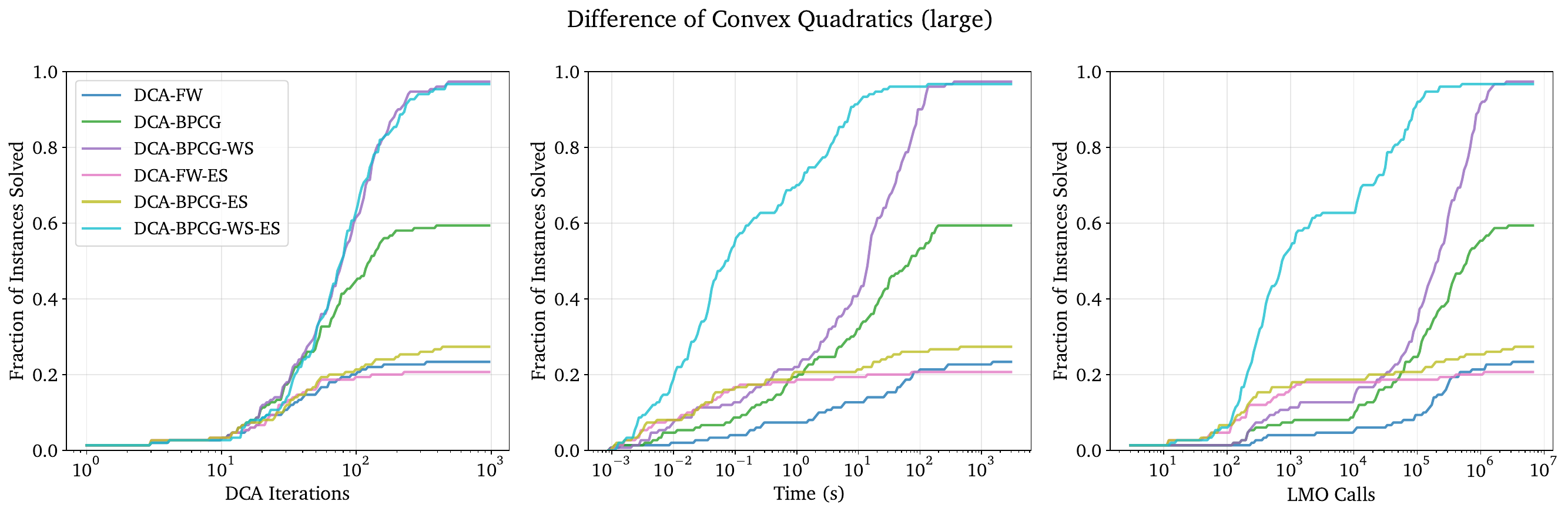}
    \caption{Performance profile over medium to large-sized instances of differences of convex quadratics from \cref{sec:diff-quadratics}.}
    \label{fig:performance-profile-simple_hard}
\end{figure}

For this straightforward experiment, we generate instances where both $f$ and $g$ are convex quadratic functions of the form
\begin{align}
    \label{eq:diff-quadratics}
    f(x) = \frac{1}{2} x^T A x + a^T x + c \quad \text{and} \quad g(x) = \frac{1}{2} x^T B x + b^T x + d,
\end{align}
where $A, B \in \mathbb{R}^{n \times n}$ are positive definite matrices (constructed as $M^T M + 0.1 I$ for random matrices $M$), $a, b \in \mathbb{R}^n$ are random vectors, and $c, d \in \mathbb{R}$ are random scalars. The feasible region is the probability simplex $P = \{x \in \mathbb{R}^n : x \geq 0, \sum_{i=1}^n x_i = 1\}$. We ran two separate tests, once with smaller to moderate instance sizes (see \cref{tab:dca_fw_statistics_simple} and \cref{fig:performance-profile-simple}), where we imposed an inner iteration limit of $10000$ iterations and a limit of $200$ DCA iterations, and once with larger instance sizes (see \cref{tab:dca_fw_statistics_simple_hard} and \cref{fig:performance-profile-simple_hard} and \cref{fig:performance-profile-simple_hard}), where we imposed an inner iteration limit of $50000$ iterations and a limit of $500$ DCA iterations. In both cases, we set the DCA gap tolerance to $10^{-6}$ and the Frank-Wolfe gap tolerance to $10^{-6}/2$.

As can be seen from the results, in both tests BPCG-WS-ES, combining warm-starting and early stopping is the most efficient outperforming any other variant by a large margin. 

\begin{table}[htbp]
\centering
\caption{Performance Statistics over small to medium-sized instances of differences of convex quadratics from \cref{sec:diff-quadratics} [Shifted Geometric Means, shift: $1.0$].}
\label{tab:dca_fw_statistics_simple}
\resizebox{\textwidth}{!}{
\begin{tabular}{lrrrrrrrrrrrrrrrrrr}
\toprule
Size & \multicolumn{3}{c}{DCA-FW} & \multicolumn{3}{c}{DCA-BPCG} & \multicolumn{3}{c}{DCA-BPCG-WS} & \multicolumn{3}{c}{DCA-FW-ES} & \multicolumn{3}{c}{DCA-BPCG-ES} & \multicolumn{3}{c}{DCA-BPCG-WS-ES} \\
 & Iter & Time & LMO & Iter & Time & LMO & Iter & Time & LMO & Iter & Time & LMO & Iter & Time & LMO & Iter & Time & LMO \\
\midrule
$n=10$ (5) & 6.36 & 0.20 & 92.98 & \textbf{4.70} & 0.05 & 19.56 & \textbf{4.70} & \textbf{2.18e-03} & 18.88 & \textbf{4.70} & 0.02 & \textbf{16.41} & \textbf{4.70} & 0.04 & \textbf{16.41} & 4.91 & 2.48e-03 & 17.17 \\
$n=20$ (5) & 75.56 & 3.84 & 6.39e+05 & \textbf{15.25} & 8.36e-03 & 886.08 & \textbf{15.25} & 2.02e-03 & 203.12 & 95.02 & 4.12 & 1.83e+05 & 95.02 & 4.47 & 1.80e+05 & 15.78 & \textbf{1.53e-03} & \textbf{79.17} \\
$n=30$ (5) & 57.96 & 1.72 & 3.26e+04 & 21.80 & 5.06e-03 & 463.72 & 22.00 & 2.18e-03 & 197.33 & 35.67 & 0.58 & 1611.83 & 22.70 & 2.25e-03 & 150.16 & \textbf{20.12} & \textbf{1.78e-03} & \textbf{116.94} \\
$n=50$ (5) & 63.54 & 1.83 & 1.02e+04 & 20.82 & 0.03 & 716.17 & \textbf{20.59} & 2.75e-03 & 186.30 & 23.52 & 5.20e-03 & 246.55 & 22.78 & 5.74e-03 & 175.84 & 21.08 & \textbf{2.32e-03} & \textbf{121.46} \\
$n=100$ (5) & 64.10 & 9.01 & 5.49e+05 & 24.93 & 0.45 & 1.56e+04 & 25.03 & 0.09 & 1173.44 & 62.98 & 5.70 & 2.82e+04 & 62.98 & 6.11 & 2.77e+04 & \textbf{22.97} & \textbf{5.29e-03} & \textbf{116.52} \\
$n=150$ (5) & 26.03 & 2.67 & 9271.30 & 22.06 & 1.42 & 3921.14 & 14.07 & 0.09 & 221.64 & 22.06 & 1.07 & 395.61 & 21.74 & 1.05 & 384.23 & \textbf{13.74} & \textbf{3.01e-03} & \textbf{58.29} \\
$n=200$ (5) & 131.48 & 44.42 & 1.26e+06 & 112.80 & 20.27 & 7.01e+05 & 58.81 & 1.35 & 1.39e+04 & 132.52 & 22.72 & 2.76e+05 & 132.52 & 23.19 & 2.61e+05 & \textbf{54.49} & \textbf{0.02} & \textbf{355.98} \\
$n=300$ (5) & 137.30 & 97.19 & 1.32e+06 & 96.82 & 19.64 & 3.70e+05 & 82.72 & 7.89 & 1.92e+05 & 139.04 & 31.12 & 2.62e+05 & 140.70 & 31.71 & 2.67e+05 & \textbf{82.72} & \textbf{0.05} & \textbf{750.10} \\
$n=400$ (5) & 200.00 & 146.02 & 1.98e+06 & 182.94 & 50.85 & 1.07e+06 & \textbf{72.22} & 6.31 & 1.16e+05 & 200.00 & 96.07 & 1.56e+06 & 200.00 & 91.76 & 1.54e+06 & 72.24 & \textbf{0.28} & \textbf{944.28} \\
$n=500$ (5) & 51.26 & 37.41 & 9.04e+04 & 38.88 & 17.77 & 5.40e+04 & \textbf{23.73} & 4.71 & 1.62e+04 & 69.20 & 20.56 & 5.52e+04 & 78.94 & 37.88 & 1.06e+05 & 37.41 & \textbf{0.40} & \textbf{865.13} \\
\bottomrule
\end{tabular}
}
\end{table}

\begin{table}[htbp]
\centering
\caption{Performance Statistics over medium to large-sized instances of differences of convex quadratics from \cref{sec:diff-quadratics} [Shifted Geometric Means, shift: $1.0$].}
\label{tab:dca_fw_statistics_simple_hard}
\resizebox{\textwidth}{!}{
\begin{tabular}{lrrrrrrrrrrrrrrrrrr}
\toprule
Size & \multicolumn{3}{c}{DCA-FW} & \multicolumn{3}{c}{DCA-BPCG} & \multicolumn{3}{c}{DCA-BPCG-WS} & \multicolumn{3}{c}{DCA-FW-ES} & \multicolumn{3}{c}{DCA-BPCG-ES} & \multicolumn{3}{c}{DCA-BPCG-WS-ES} \\
 & Iter & Time & LMO & Iter & Time & LMO & Iter & Time & LMO & Iter & Time & LMO & Iter & Time & LMO & Iter & Time & LMO \\
\midrule
$n=100$ (10) & 231.25 & 23.32 & 9.58e+05 & 47.31 & 0.47 & 2.03e+04 & 47.17 & 0.11 & 2362.48 & 124.64 & 7.56 & 7.04e+04 & 169.42 & 12.94 & 1.74e+05 & \textbf{37.85} & \textbf{6.84e-03} & \textbf{230.11} \\
$n=200$ (10) & 295.57 & 85.48 & 2.54e+06 & 108.16 & 8.94 & 2.81e+05 & \textbf{76.29} & 2.16 & 5.71e+04 & 294.68 & 61.78 & 1.06e+06 & 297.27 & 68.48 & 1.13e+06 & 84.33 & \textbf{0.08} & \textbf{981.80} \\
$n=300$ (10) & 360.18 & 188.71 & 3.51e+06 & 144.27 & 27.36 & 6.37e+05 & \textbf{73.15} & 5.60 & 1.42e+05 & 365.52 & 138.48 & 1.84e+06 & 500.00 & 240.35 & 4.61e+06 & 81.69 & \textbf{0.25} & \textbf{1985.46} \\
$n=400$ (10) & 405.53 & 218.76 & 1.85e+06 & 337.29 & 116.75 & 1.38e+06 & 95.30 & 9.11 & 7.60e+04 & 390.92 & 175.82 & 1.65e+06 & 390.92 & 176.12 & 1.67e+06 & \textbf{94.94} & \textbf{0.71} & \textbf{3384.47} \\
$n=500$ (10) & 398.44 & 402.08 & 3.92e+06 & 272.77 & 142.85 & 2.23e+06 & 101.70 & 20.20 & 3.62e+05 & 500.00 & 403.29 & 4.47e+06 & 406.82 & 224.01 & 2.54e+06 & \textbf{101.19} & \textbf{1.06} & \textbf{1.17e+04} \\
$n=600$ (10) & 400.69 & 408.89 & 2.33e+06 & 217.09 & 135.01 & 1.70e+06 & \textbf{121.52} & 30.45 & 3.78e+05 & 399.20 & 306.11 & 1.86e+06 & 402.62 & 222.93 & 1.39e+06 & 135.39 & \textbf{1.35} & \textbf{7345.49} \\
$n=700$ (10) & 421.90 & 298.50 & 1.92e+06 & 298.63 & 109.60 & 9.26e+05 & 92.27 & 12.81 & 9.15e+04 & 399.95 & 260.05 & 1.74e+06 & 306.37 & 127.02 & 9.67e+05 & \textbf{80.22} & \textbf{1.14} & \textbf{3409.07} \\
$n=800$ (10) & 365.52 & 579.95 & 3.42e+06 & 356.54 & 346.69 & 3.09e+06 & 128.18 & 40.42 & 3.97e+05 & 368.72 & 251.60 & 1.38e+06 & 360.66 & 209.15 & 1.35e+06 & \textbf{125.44} & \textbf{1.66} & \textbf{5473.48} \\
$n=900$ (10) & 221.82 & 331.17 & 8.10e+05 & 159.99 & 191.99 & 6.21e+05 & 88.35 & 46.95 & 1.90e+05 & 200.90 & 121.79 & 2.68e+05 & 170.29 & 87.36 & 2.10e+05 & \textbf{86.70} & \textbf{2.21} & \textbf{5270.98} \\
$n=1000$ (10) & 236.52 & 239.21 & 4.59e+05 & 151.12 & 90.62 & 2.32e+05 & 88.63 & 27.89 & 7.70e+04 & 228.92 & 204.13 & 3.94e+05 & 197.63 & 129.56 & 2.79e+05 & \textbf{88.62} & \textbf{1.92} & \textbf{3084.11} \\
\bottomrule
\end{tabular}
}
\end{table}

\subsection{Hard DC function example}
\label{sec:hard-dc}

As a more challenging example we consider a DC optimization problem that combines quadratic, exponential, and logistic components. Specifically, we minimize $\dcprob{x} = f(x) - g(x)$ where 
$$f(x) = \frac{1}{2} x^T A x + a^T x + \frac{1}{n} \exp\left(\frac{1}{n} \sum_{i=1}^n c_i x_i\right),$$
and 
$$g(x) = \frac{1}{2} x^T B x + b^T x + 0.1 \sum_{i=1}^n \log(1 + \exp(d_i x_i)),$$
over the $k$-sparse polytope $C_{\tau, k} = B_1(k \tau) \cap B_\infty(\tau) =  \{x \in \mathbb{R}^n : \|x\|_1 \leq k \tau, \|x\|_\infty \leq \tau\}$ with $\tau = 10$ and $k = 10$, which is the intersection of the $\ell_1$-ball of radius $k \tau$ and the $\ell_\infty$-ball of radius $\tau$. This example is particularly challenging because the exponential term in $f(x)$ can lead to rapid growth, while the logistic terms in $g(x)$ introduce additional nonlinearity. The $k$-sparse polytope further complicates the optimization landscape. The matrices $A$, $B$ and the vectors $a, b, c, d \in \mathbb{R}^n$ are constructed as in \cref{sec:diff-quadratics}. 

\begin{figure}[htbp]
    \centering
        \includegraphics[width=0.99\textwidth]{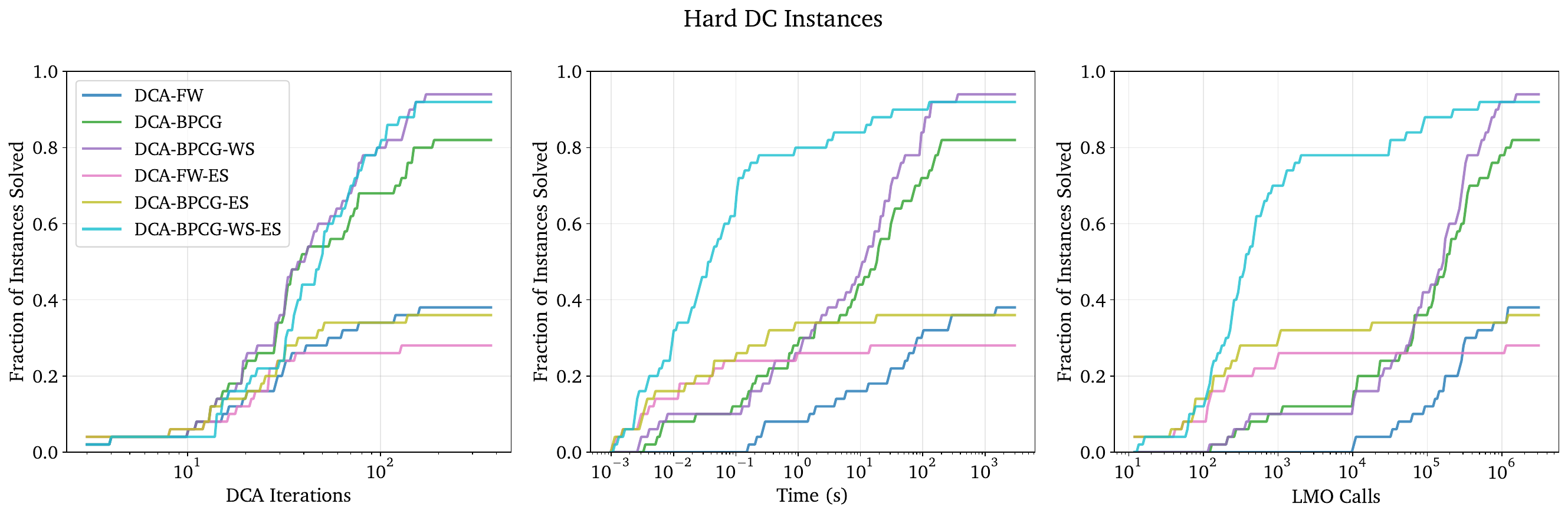}
    \caption{Performance profile over hard DC instances from \cref{sec:hard-dc}.}
    \label{fig:performance-profile-hard_dc}
\end{figure}

Results are reported in \cref{tab:dca_fw_statistics_hard_dc} and \cref{fig:performance-profile-hard_dc} and match the observations from \cref{sec:diff-quadratics}.

\begin{table}[htbp]
\centering
\caption{Performance Statistics for hard DC instances from \cref{sec:hard-dc} [Shifted Geometric Means, shift: $1.0$].}
\label{tab:dca_fw_statistics_hard_dc}
\resizebox{\textwidth}{!}{
\begin{tabular}{lrrrrrrrrrrrrrrrrrr}
\toprule
Size & \multicolumn{3}{c}{DCA-FW} & \multicolumn{3}{c}{DCA-BPCG} & \multicolumn{3}{c}{DCA-BPCG-WS} & \multicolumn{3}{c}{DCA-FW-ES} & \multicolumn{3}{c}{DCA-BPCG-ES} & \multicolumn{3}{c}{DCA-BPCG-WS-ES} \\
 & Iter & Time & LMO & Iter & Time & LMO & Iter & Time & LMO & Iter & Time & LMO & Iter & Time & LMO & Iter & Time & LMO \\
\midrule
$n=10$ (5) & 88.42 & 4.53 & 4.03e+05 & \textbf{16.21} & 0.06 & 557.81 & \textbf{16.21} & 0.07 & 675.40 & 76.80 & 2.03 & 1.09e+04 & 52.38 & 2.28 & 1.35e+04 & 18.21 & \textbf{3.60e-03} & \textbf{101.01} \\
$n=20$ (5) & 54.03 & 6.94 & 3.88e+05 & \textbf{16.91} & 0.34 & 1.27e+04 & 17.30 & 0.35 & 1.17e+04 & 53.47 & 5.72 & 3.33e+04 & 53.21 & 6.03 & 2.29e+04 & 17.26 & \textbf{4.14e-03} & \textbf{102.78} \\
$n=30$ (5) & 66.14 & 9.44 & 3.41e+05 & 20.41 & 0.63 & 1.76e+04 & 20.54 & 0.89 & 3.61e+04 & 23.61 & 1.07 & 692.70 & 36.84 & 3.58 & 3300.69 & \textbf{19.17} & \textbf{0.01} & \textbf{203.28} \\
$n=50$ (5) & 71.10 & 30.98 & 6.57e+05 & 57.88 & 8.69 & 2.57e+05 & \textbf{25.00} & 2.05 & 6.55e+04 & 107.00 & 36.66 & 2.06e+05 & 107.00 & 29.27 & 1.88e+05 & 43.39 & \textbf{0.15} & \textbf{894.76} \\
$n=100$ (5) & 126.88 & 137.60 & 1.22e+06 & 71.12 & 21.09 & 3.36e+05 & 53.00 & 14.78 & 2.37e+05 & 200.00 & 166.19 & 1.60e+06 & 136.40 & 79.22 & 2.86e+05 & \textbf{51.84} & \textbf{0.40} & \textbf{1062.69} \\
$n=150$ (5) & 140.70 & 237.07 & 1.39e+06 & 70.91 & 30.58 & 3.87e+05 & \textbf{48.67} & 13.28 & 1.68e+05 & 133.53 & 95.10 & 2.80e+05 & 116.97 & 77.18 & 2.22e+05 & 54.40 & \textbf{0.04} & \textbf{312.59} \\
$n=200$ (5) & 50.51 & 117.98 & 4.58e+05 & \textbf{32.48} & 15.49 & 1.49e+05 & 34.58 & 14.64 & 1.38e+05 & 82.93 & 37.51 & 4.22e+04 & 63.52 & 11.91 & 8853.26 & 43.61 & \textbf{0.36} & \textbf{744.93} \\
$n=300$ (5) & 135.47 & 583.05 & 1.33e+06 & 114.63 & 136.54 & 9.50e+05 & \textbf{98.29} & 94.91 & 6.53e+05 & 200.00 & 379.09 & 8.90e+05 & 132.52 & 120.06 & 1.45e+05 & 121.35 & \textbf{2.73} & \textbf{6649.18} \\
$n=400$ (5) & 138.18 & 1042.07 & 1.36e+06 & 108.40 & 111.16 & 5.99e+05 & 108.57 & 98.88 & 5.18e+05 & 200.00 & 908.32 & 1.30e+06 & 149.94 & 160.82 & 1.85e+05 & \textbf{65.05} & \textbf{0.91} & \textbf{1146.13} \\
$n=500$ (5) & 130.43 & 1170.53 & 1.20e+06 & 101.83 & 107.13 & 3.83e+05 & \textbf{99.72} & 99.51 & 3.77e+05 & 142.28 & 305.00 & 3.06e+05 & 139.04 & 162.55 & 1.59e+05 & 152.95 & \textbf{20.07} & \textbf{3.96e+04} \\
\bottomrule
\end{tabular}
}
\end{table}

\subsection{Quadratic Assignment Problem}
\label{sec:qap}

In this test we solve a DC relaxation of the Quadratic Assignment Problem (QAP). The QAP aims to find an optimal alignment between two matrices $A, B \in \mathbb{R}^{n\times n}$, and it is well known that this is equivalent to minimizing a quadratic objective function over the set of permutation matrices
\begin{align}
\label{eqn:QAP}
   \underset{X \in \mathbb{R}^{n \times n}}{\min} \innp{A^\top X, X B} \quad
   \mathrm{s.t.}  \quad  X \in \{0,1\}^{n \times n}, ~ X1_n = X^\top1_n = 1_n,
\end{align}
where $1_n$ denotes the $n$-dimensional vector of ones. As this combinatorial optimization problem is NP-hard, we want to solve it over the linear relaxation of the constraints, where we drop the integrality required:
\begin{align}
\label{eqn:QAP-relaxed}
   \underset{X \in \mathbb{R}^{n \times n}}{\min} ~~  \innp{A^\top X, X B} \quad
   \mathrm{s.t.}  \quad  X \in [0,1]^{n \times n}, ~ X1_n = X^\top1_n = 1_n,
\end{align}
i.e., we want to minimize a nonconvex quadratic over the Birkhoff polytope $\mathcal{B}_n = \{X \in \mathbb{R}^{n \times n} : X 1_n = 1_n, X^\top 1_n = 1_n, X \geq 0\}$. This solution can then often be used for rounding to obtain an actual feasible assignment of the original problem; see, e.g., \citet{vogelstein2015fast}, see also \citet{maskanRevisitingFrankWolfeStructured2025} for an in-depth discussion. 

In order to apply DCA, we need to reformulate the QAP and in fact decompose the objective $\innp{A^\top X, X B}$ into a difference of convex functions. Several decompositions are possible and we opted for the approach in \citet{maskanRevisitingFrankWolfeStructured2025} based on their tests and for comparability:
$$f(X) = \frac{1}{4}\|A^\top X + XB\|_F^2,$$
and
$$g(X) = \frac{1}{4}\|A^\top X - XB\|_F^2,$$
where $\|\cdot\|_F$ denotes the Frobenius norm. This gives us $\innp{A^\top X, X B} = f(X) - g(X)$ since
$$f(X) - g(X) = \frac{1}{4}\left(\|A^\top X + XB\|_F^2 - \|A^\top X - XB\|_F^2\right) = \innp{A^\top X, X B}.$$
Note that both $f$ and $g$ are convex functions as compositions of convex quadratic forms with linear transformations.

While projection onto the Birkhoff polytope is computationally expensive, the \lmo can be implemented efficiently using the Hungarian algorithm, which is readily available in the \frankwolfejl package. We evaluate the algorithms on instances from QAPLIB \citep{burkard1997qaplib}, a well-established benchmark library for quadratic assignment problems\footnote{available at \url{https://coral.ise.lehigh.edu/data-sets/qaplib/qaplib-problem-instances-and-solutions/}}. QAPLIB contains a wide range of QAP instances that have been used extensively in the literature to benchmark QAP algorithms. From the 136 instances in QAPLIB, we successfully processed 97 instances (71.3\%) for our experiments. The remaining 39 instances (28.7\%) could not be processed due to various data format inconsistencies or parsing issues; for details see \cref{sec:qap-instances}.

For our experiments, we used a DCA iteration limit of $500$ and an inner iteration limit for the Frank-Wolfe iterations of $50000$. The instance characteristics can be found in \cref{tab:dca_fw_statistics_qap}. Note that the actual problem dimension is quadratic in $n$. Due to their considerable size and the substantially higher \lmo complexity (as we have to now run the Hungarian algorithm compared to e.g., sorting before), these instances are much harder than the previous examples, also often coming with objective values that are orders of magnitude larger, so that the additive tolerance is relatively tighter. Thus, the target DCA tolerance of $10^{-6}$ could not always be achieved, even for the fastest algorithm configuration.

\begin{figure}[htbp]
    \centering
        \includegraphics[width=0.99\textwidth]{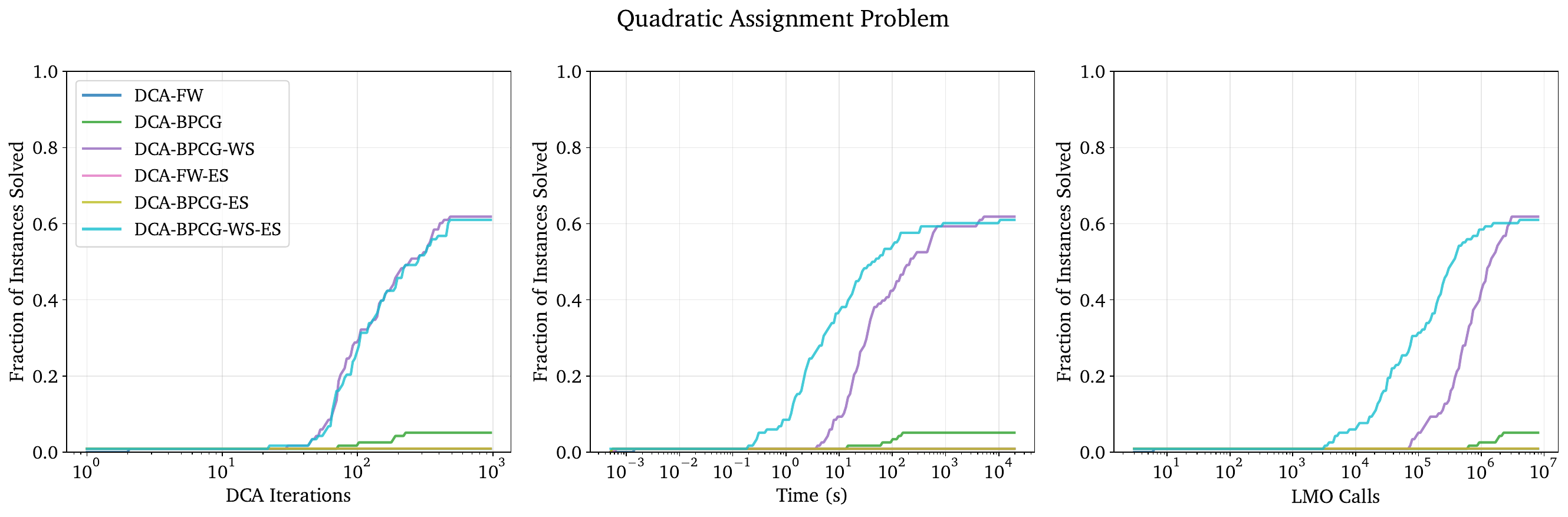}
    \caption{Performance profile over Quadratic Assignment Problem instances from \cref{sec:qap}.}
    \label{fig:performance-profile-qap}
\end{figure}

\begin{figure}[htbp]
    \centering
        \includegraphics[width=0.99\textwidth]{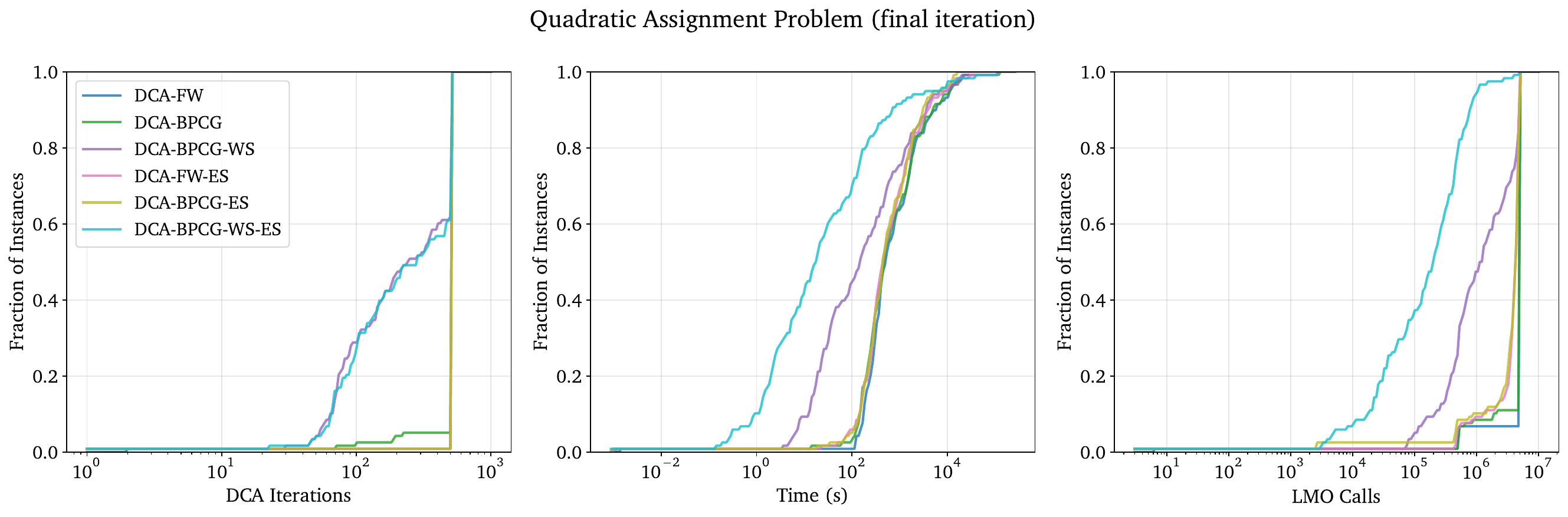}
    \caption{Modified performance profile over Quadratic Assignment Problem instances from \cref{sec:qap}, where we consider all instances to be solved in their last iteration.}
    \label{fig:performance-profile-qap-final}
\end{figure}

Given these computational challenges, we present two complementary performance profiles for these experiments: the classical performance profile in \cref{fig:performance-profile-qap}, and a modified performance profile in \cref{fig:performance-profile-qap-final} where we consider all instances to be ``solved'' at their final iteration. This second profile allows us to compare the computational cost---both in terms of runtime and \lmo calls---required to reach the DCA iteration limit across different algorithm variants.

Also note that in \cref{tab:dca_fw_statistics_qap}, the results seem to be a little bit more noisy, which in fact is due to often having only one QAP instance in a given size bracket. Due to nonconvexity, the optimization paths can differ signficantly, so that sometimes the slower algorithms can be faster as they converged to a different stationary point.

\begin{table}[htbp]
\centering
\caption{Performance Statistics Quadratic Assignment Problem instances from \cref{sec:qap} [Shifted Geometric Means, shift: $1.0$].}
\label{tab:dca_fw_statistics_qap}
\resizebox{\textwidth}{!}{
\begin{tabular}{lrrrrrrrrrrrrrrrrrr}
\toprule
Size & \multicolumn{3}{c}{DCA-FW} & \multicolumn{3}{c}{DCA-BPCG} & \multicolumn{3}{c}{DCA-BPCG-WS} & \multicolumn{3}{c}{DCA-FW-ES} & \multicolumn{3}{c}{DCA-BPCG-ES} & \multicolumn{3}{c}{DCA-BPCG-WS-ES} \\
 & Iter & Time & LMO & Iter & Time & LMO & Iter & Time & LMO & Iter & Time & LMO & Iter & Time & LMO & Iter & Time & LMO \\
\midrule
$n=10$ (4) & 500.00 & 120.53 & 5.00e+06 & 307.47 & 70.15 & 2.96e+06 & 177.98 & 31.30 & 1.43e+06 & 500.00 & 88.75 & 3.52e+06 & 500.00 & 105.81 & 3.76e+06 & \textbf{168.84} & \textbf{1.85} & \textbf{7.42e+04} \\
$n=12$ (18) & 500.00 & 193.50 & 5.00e+06 & 500.00 & 182.21 & 4.96e+06 & 164.44 & 36.37 & 9.97e+05 & 500.00 & 149.91 & 3.71e+06 & 500.00 & 159.10 & 3.71e+06 & \textbf{161.41} & \textbf{3.83} & \textbf{8.45e+04} \\
$n=14$ (4) & 500.00 & 286.42 & 5.00e+06 & 500.00 & 258.41 & 4.95e+06 & \textbf{96.73} & 12.54 & 2.48e+05 & 500.00 & 270.04 & 4.57e+06 & 500.00 & 287.84 & 4.63e+06 & 101.40 & \textbf{2.79} & \textbf{4.82e+04} \\
$n=15$ (16) & 500.00 & 315.63 & 5.00e+06 & 474.88 & 280.30 & 4.71e+06 & \textbf{201.66} & 79.17 & 1.37e+06 & 500.00 & 246.40 & 3.75e+06 & 500.00 & 230.92 & 3.30e+06 & 241.32 & \textbf{12.58} & \textbf{1.51e+05} \\
$n=16$ (13) & 336.96 & 243.49 & 1.77e+06 & 280.60 & 195.65 & 1.44e+06 & 58.17 & 12.45 & 9.34e+04 & 326.58 & 223.25 & 1.53e+06 & 326.58 & 242.75 & 1.55e+06 & \textbf{52.80} & \textbf{1.25} & \textbf{8659.07} \\
$n=17$ (2) & 500.00 & 402.47 & 5.00e+06 & 500.00 & 334.74 & 4.89e+06 & \textbf{80.39} & 30.59 & 4.50e+05 & 500.00 & 376.93 & 4.64e+06 & 500.00 & 394.03 & 4.56e+06 & 88.36 & \textbf{4.96} & \textbf{6.14e+04} \\
$n=18$ (4) & 500.00 & 457.20 & 5.00e+06 & 334.71 & 267.48 & 3.22e+06 & \textbf{141.19} & 41.64 & 5.01e+05 & 500.00 & 444.26 & 4.61e+06 & 500.00 & 460.44 & 4.56e+06 & 216.88 & \textbf{10.29} & \textbf{1.16e+05} \\
$n=20$ (8) & 500.00 & 629.25 & 5.00e+06 & 500.00 & 574.35 & 4.97e+06 & 228.20 & 148.40 & 1.35e+06 & 500.00 & 471.76 & 3.66e+06 & 500.00 & 543.62 & 3.96e+06 & \textbf{219.57} & \textbf{11.92} & \textbf{9.83e+04} \\
$n=21$ (1) & 500.00 & 665.66 & 5.00e+06 & 500.00 & 599.23 & 4.92e+06 & \textbf{300.00} & 109.76 & 9.54e+05 & 500.00 & 539.76 & 3.77e+06 & 500.00 & 545.62 & 3.61e+06 & 455.00 & \textbf{49.25} & \textbf{4.09e+05} \\
$n=22$ (3) & 500.00 & 695.69 & 5.00e+06 & 500.00 & 620.95 & 4.95e+06 & \textbf{298.58} & 179.39 & 1.53e+06 & 500.00 & 570.54 & 4.02e+06 & 500.00 & 612.17 & 4.19e+06 & 411.68 & \textbf{51.99} & \textbf{4.19e+05} \\
$n=24$ (1) & 500.00 & 865.14 & 5.00e+06 & 500.00 & 758.07 & 4.93e+06 & 432.00 & 180.65 & 1.24e+06 & 500.00 & 704.49 & 3.95e+06 & 500.00 & 721.46 & 3.86e+06 & \textbf{101.00} & \textbf{7.98} & \textbf{4.99e+04} \\
$n=25$ (4) & 500.00 & 1096.12 & 5.00e+06 & 500.00 & 957.49 & 4.93e+06 & \textbf{350.22} & 458.78 & 2.55e+06 & 500.00 & 832.00 & 3.76e+06 & 500.00 & 939.65 & 4.05e+06 & 370.15 & \textbf{76.56} & \textbf{4.22e+05} \\
$n=26$ (6) & \textbf{500.00} & 1472.84 & 5.00e+06 & \textbf{500.00} & 1522.93 & 5.00e+06 & \textbf{500.00} & 1262.66 & 4.93e+06 & \textbf{500.00} & 1097.81 & 3.55e+06 & \textbf{500.00} & 1193.52 & 3.63e+06 & \textbf{500.00} & \textbf{163.96} & \textbf{6.14e+05} \\
$n=27$ (1) & 500.00 & 1109.44 & 5.00e+06 & 500.00 & 983.31 & 4.96e+06 & 202.00 & 205.24 & 1.06e+06 & 500.00 & 927.43 & 3.87e+06 & 500.00 & 1106.67 & 4.64e+06 & \textbf{196.00} & \textbf{39.40} & \textbf{1.88e+05} \\
$n=28$ (1) & 500.00 & 1267.68 & 5.00e+06 & 500.00 & 1076.75 & 4.95e+06 & \textbf{236.00} & 137.59 & 6.31e+05 & 500.00 & 853.08 & 3.33e+06 & 500.00 & 749.71 & 2.75e+06 & 346.00 & \textbf{20.35} & \textbf{9.48e+04} \\
$n=30$ (3) & 500.00 & 1931.15 & 5.00e+06 & 500.00 & 1788.34 & 4.96e+06 & 463.39 & 1144.63 & 3.56e+06 & 500.00 & 1485.99 & 3.91e+06 & 500.00 & 1518.37 & 3.86e+06 & \textbf{377.48} & \textbf{107.08} & \textbf{3.33e+05} \\
$n=32$ (8) & 500.00 & 1753.54 & 5.00e+06 & 500.00 & 1683.56 & 5.00e+06 & 149.63 & 459.75 & 1.44e+06 & 500.00 & 1606.51 & 4.43e+06 & 500.00 & 1649.78 & 4.32e+06 & \textbf{129.19} & \textbf{91.05} & \textbf{2.63e+05} \\
$n=35$ (2) & 500.00 & 2507.13 & 5.00e+06 & 500.00 & 2208.98 & 4.96e+06 & \textbf{324.90} & 1310.80 & 2.98e+06 & 500.00 & 1680.52 & 3.47e+06 & 500.00 & 2370.69 & 4.43e+06 & 500.00 & \textbf{81.36} & \textbf{1.73e+05} \\
$n=40$ (3) & 500.00 & 3294.19 & 5.00e+06 & 500.00 & 2908.98 & 4.97e+06 & 365.42 & 1622.86 & 2.89e+06 & 500.00 & 2707.80 & 4.01e+06 & 500.00 & 2625.10 & 3.96e+06 & \textbf{341.45} & \textbf{206.14} & \textbf{3.72e+05} \\
$n=50$ (3) & \textbf{500.00} & 5879.99 & 5.00e+06 & \textbf{500.00} & 5006.42 & 4.96e+06 & \textbf{500.00} & 3038.16 & 3.15e+06 & \textbf{500.00} & 4199.93 & 3.73e+06 & \textbf{500.00} & 3423.85 & 2.82e+06 & \textbf{500.00} & \textbf{341.62} & \textbf{3.50e+05} \\
$n=60$ (2) & 500.00 & 8838.29 & 5.00e+06 & 500.00 & 7606.62 & 4.96e+06 & 500.00 & 6400.82 & 4.31e+06 & 500.00 & 6513.19 & 4.01e+06 & 500.00 & 6672.89 & 3.61e+06 & \textbf{417.15} & \textbf{1106.33} & \textbf{7.34e+05} \\
$n=64$ (2) & 500.00 & 1.48e+04 & 5.00e+06 & 500.00 & 1.47e+04 & 5.00e+06 & \textbf{334.75} & 8705.74 & 3.35e+06 & 500.00 & 1.79e+04 & 4.64e+06 & 500.00 & \textbf{506.80} & \textbf{1.08e+05} & 488.36 & 1.93e+04 & 4.36e+06 \\
$n=80$ (2) & 500.00 & 5185.01 & 1.58e+06 & 500.00 & 4635.28 & 1.57e+06 & \textbf{415.35} & 3062.99 & 1.13e+06 & 500.00 & 1968.18 & 6.12e+05 & 500.00 & 2826.78 & 8.21e+05 & 500.00 & \textbf{1467.19} & \textbf{5.74e+05} \\
$n=100$ (3) & \textbf{500.00} & 6667.74 & 1.08e+06 & \textbf{500.00} & 6153.12 & 1.08e+06 & \textbf{500.00} & 5643.66 & 1.03e+06 & \textbf{500.00} & 3912.33 & 5.77e+05 & \textbf{500.00} & 3559.44 & 5.60e+05 & \textbf{500.00} & \textbf{596.97} & \textbf{1.13e+05} \\
$n=128$ (1) & 500.00 & 7165.42 & 5.02e+05 & 500.00 & 8143.55 & 5.02e+05 & \textbf{400.00} & \textbf{5057.93} & \textbf{4.01e+05} & 500.00 & 8224.50 & 4.88e+05 & 500.00 & 9925.79 & 4.85e+05 & 500.00 & 7541.98 & 4.81e+05 \\
$n=150$ (2) & \textbf{500.00} & 1.35e+04 & 5.02e+05 & \textbf{500.00} & 1.46e+04 & 5.02e+05 & \textbf{500.00} & 1.34e+04 & 5.02e+05 & \textbf{500.00} & 1.31e+04 & 4.79e+05 & \textbf{500.00} & 1.31e+04 & 4.55e+05 & \textbf{500.00} & \textbf{1.24e+04} & \textbf{4.52e+05} \\
$n=256$ (1) & \textbf{500.00} & 1.16e+05 & 5.02e+05 & \textbf{500.00} & 1.30e+05 & 5.02e+05 & \textbf{500.00} & 1.17e+05 & 5.02e+05 & \textbf{500.00} & 1.16e+05 & 4.96e+05 & \textbf{500.00} & \textbf{612.81} & \textbf{2500.00} & \textbf{500.00} & 1.11e+05 & 4.59e+05 \\
\bottomrule
\end{tabular}
}
\end{table}

\section*{Acknowledgements}

Research reported in this paper was partially supported through the Research Campus Modal funded by the German Federal Ministry of Education and Research (fund numbers 05M14ZAM,05M20ZBM) and the Deutsche Forschungsgemeinschaft (DFG) through the DFG Cluster of Excellence MATH+. 

\bibliography{pubPokutta,refs}
\bibliographystyle{icml2021}

\clearpage

\appendix

\section{QAP instances}
\label{sec:qap-instances}

We evaluate our algorithms on instances from QAPLIB \citep{burkard1997qaplib}, a well-established benchmark library for quadratic assignment problems available at \url{https://coral.ise.lehigh.edu/data-sets/qaplib/qaplib-problem-instances-and-solutions/}. QAPLIB contains a wide range of QAP instances that have been used extensively in the literature to benchmark QAP algorithms.

From the 136 instances in QAPLIB, we successfully processed 97 instances (71.3\%) for our experiments. The remaining 39 instances (28.7\%) could not be processed due to various data format inconsistencies or parsing issues. The invalid instances include: \texttt{bur26c}, \texttt{bur26d}, \texttt{els19}, \texttt{kra30a}, \texttt{kra30b}, \texttt{lipa20a}, \texttt{lipa20b}, \texttt{lipa30a}, \texttt{lipa30b}, \texttt{lipa40a}, \texttt{lipa40b}, \texttt{lipa50a}, \texttt{lipa50b}, \texttt{lipa60a}, \texttt{lipa60b}, \texttt{lipa70a}, \texttt{lipa70b}, \texttt{lipa80a}, \texttt{lipa80b}, \texttt{lipa90a}, \texttt{lipa90b}, \texttt{nug30}, \texttt{scr20}, \texttt{sko100a}, \texttt{sko100b}, \texttt{sko100c}, \texttt{sko100d}, \texttt{sko100e}, \texttt{sko100f}, \texttt{sko42}, \texttt{sko49}, \texttt{sko56}, \texttt{sko64}, \texttt{sko72}, \texttt{sko81}, \texttt{sko90}, \texttt{ste36a}, \texttt{ste36b}, and \texttt{ste36c}.

The 97 valid instances include diverse problem families spanning multiple size ranges:
\begin{itemize}
\item \textbf{Small instances} ($n \leq 15$): 21 instances, including \texttt{chr12a}-\texttt{chr15c}, \texttt{esc16a}-\texttt{esc16j}, \texttt{had12}-\texttt{had14}, \texttt{nug12}-\texttt{nug15}, \texttt{rou12}-\texttt{rou15}, \texttt{scr12}-\texttt{scr15}, \texttt{tai10a}-\texttt{tai15b}.
\item \textbf{Medium instances} ($16 \leq n \leq 30$): 47 instances, including the \texttt{bur26}, \texttt{chr}, \texttt{esc32}, \texttt{had16}-\texttt{had20}, \texttt{nug16a}-\texttt{nug28}, \texttt{rou20}, \texttt{tai17a}-\texttt{tai30b}, and \texttt{tho30} families.
\item \textbf{Large instances} ($n \geq 31$): 29 instances, including \texttt{esc64a}, \texttt{esc128}, \texttt{kra32}, \texttt{tai35a}-\texttt{tai256c}, \texttt{tho40}, \texttt{tho150}, \texttt{wil50}, and \texttt{wil100}.
\end{itemize}

The complete list of 97 valid instances processed in our experiments is: \texttt{bur26a}, \texttt{bur26b}, \texttt{bur26e}, \texttt{bur26f}, \texttt{bur26g}, \texttt{bur26h}, \texttt{chr12a}, \texttt{chr12b}, \texttt{chr12c}, \texttt{chr15a}, \texttt{chr15b}, \texttt{chr15c}, \texttt{chr18a}, \texttt{chr18b}, \texttt{chr20a}, \texttt{chr20b}, \texttt{chr20c}, \texttt{chr22a}, \texttt{chr22b}, \texttt{chr25a}, \texttt{esc128}, \texttt{esc16a}, \texttt{esc16b}, \texttt{esc16c}, \texttt{esc16d}, \texttt{esc16e}, \texttt{esc16f}, \texttt{esc16g}, \texttt{esc16h}, \texttt{esc16i}, \texttt{esc16j}, \texttt{esc32a}, \texttt{esc32b}, \texttt{esc32c}, \texttt{esc32d}, \texttt{esc32e}, \texttt{esc32g}, \texttt{esc32h}, \texttt{esc64a}, \texttt{had12}, \texttt{had14}, \texttt{had16}, \texttt{had18}, \texttt{had20}, \texttt{kra32}, \texttt{nug12}, \texttt{nug14}, \texttt{nug15}, \texttt{nug16a}, \texttt{nug16b}, \texttt{nug17}, \texttt{nug18}, \texttt{nug20}, \texttt{nug21}, \texttt{nug22}, \texttt{nug24}, \texttt{nug25}, \texttt{nug27}, \texttt{nug28}, \texttt{rou12}, \texttt{rou15}, \texttt{rou20}, \texttt{scr12}, \texttt{scr15}, \texttt{tai100a}, \texttt{tai100b}, \texttt{tai10a}, \texttt{tai10b}, \texttt{tai12a}, \texttt{tai12b}, \texttt{tai150b}, \texttt{tai15a}, \texttt{tai15b}, \texttt{tai17a}, \texttt{tai20a}, \texttt{tai20b}, \texttt{tai256c}, \texttt{tai25a}, \texttt{tai25b}, \texttt{tai30a}, \texttt{tai30b}, \texttt{tai35a}, \texttt{tai35b}, \texttt{tai40a}, \texttt{tai40b}, \texttt{tai50a}, \texttt{tai50b}, \texttt{tai60a}, \texttt{tai60b}, \texttt{tai64c}, \texttt{tai80a}, \texttt{tai80b}, \texttt{tho150}, \texttt{tho30}, \texttt{tho40}, \texttt{wil100}, \texttt{wil50}.

\section{Boosted DCA}
\label{sec:boosted-dc}

\begin{figure}[htbp]
    \centering
        \includegraphics[width=0.99\textwidth]{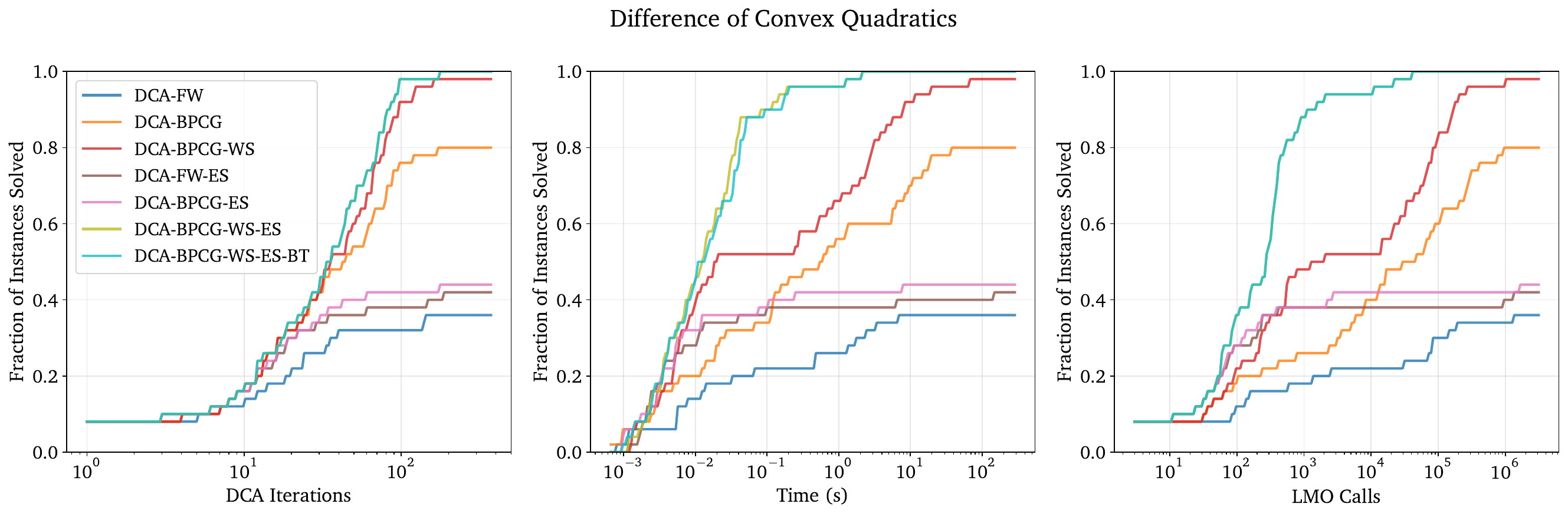}
    \caption{Performance profile over difference of convex quadratic functions as initial test of the effect of boosted DCA. The boosted DCA variant is denoted by DCA-BPCG-WS-ES-BT and we can see that its performance is basically identical to the unboosted DCA-BPCG-WS-ES variant.}
    \label{fig:performance-profile-boosted}
\end{figure}

In this section we present the preliminary results of the boosted DCA variant. We implemented a simple two level line search via refined grids to compute the step size; relatively exact but not the most efficient. The results are presented in \cref{fig:performance-profile-boosted}. We can see that the boosted DCA variant is basically identical to the unboosted DCA-BPCG-WS-ES variant. Neither is there a significant improvement in DCA iterations nor is it much slower in time due to the additional line search. Note however, that these are preliminary results and an in-depth analysis is left for future work.

\section{Additional Numerical Experiments}
\label{sec:additional-numerical-experiments}

In the following we present some selected runs to complement the summary statistics from before with some more details.

\subsection{Difference of Convex Quadratic Functions}

\begin{figure}[htbp]
    \centering
        \includegraphics[width=0.9\textwidth]{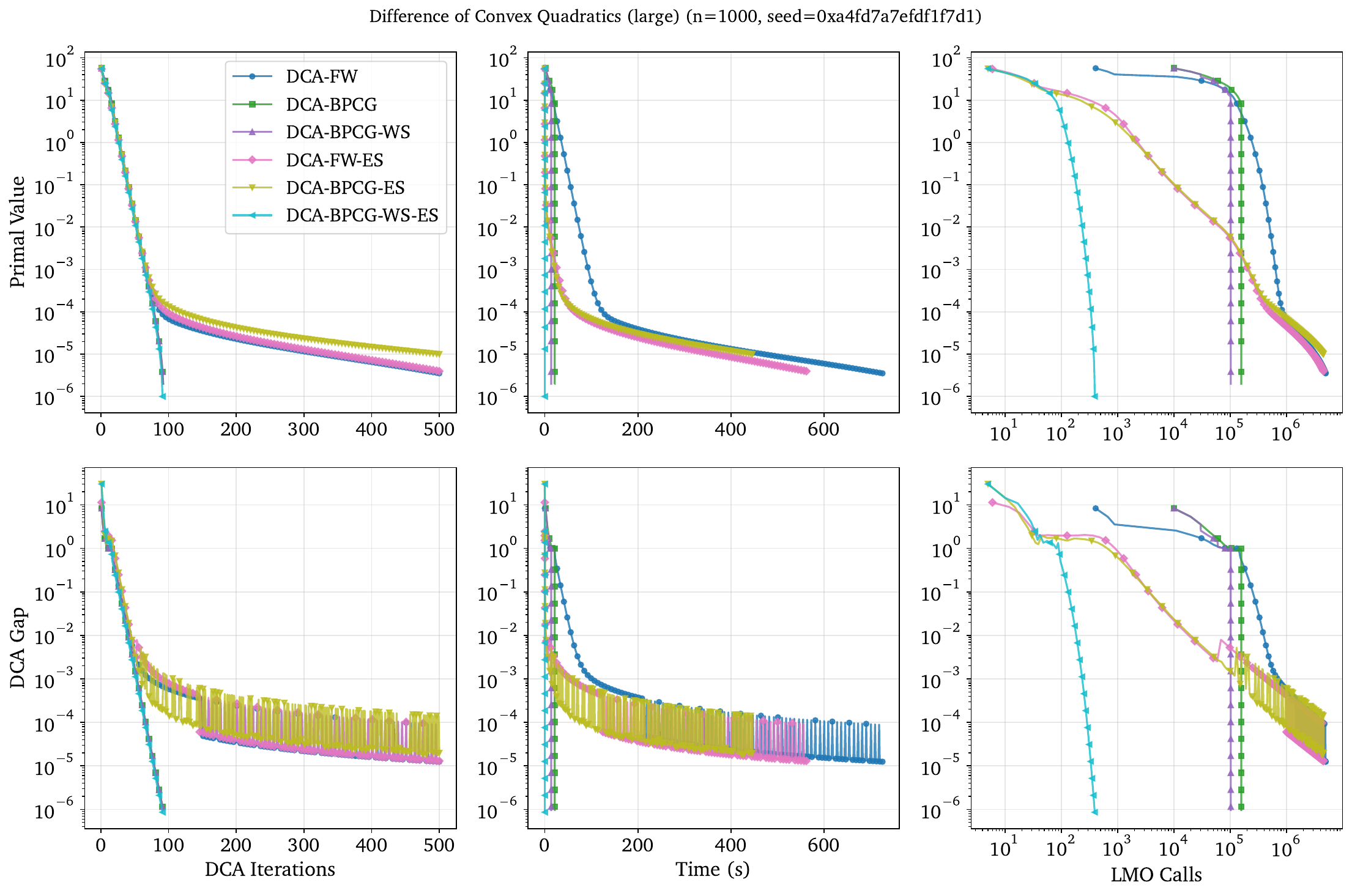}
    \caption{Run of the different variants on a Difference of Convex Quadratics example with seed $a4fd7a7efdf1f7d1$ with size $n = 1000$ from \cref{sec:diff-quadratics}.}
    \label{fig:diff-quadratics-1000-1}
\end{figure}

\begin{figure}[htbp]
    \centering
        \includegraphics[width=0.9\textwidth]{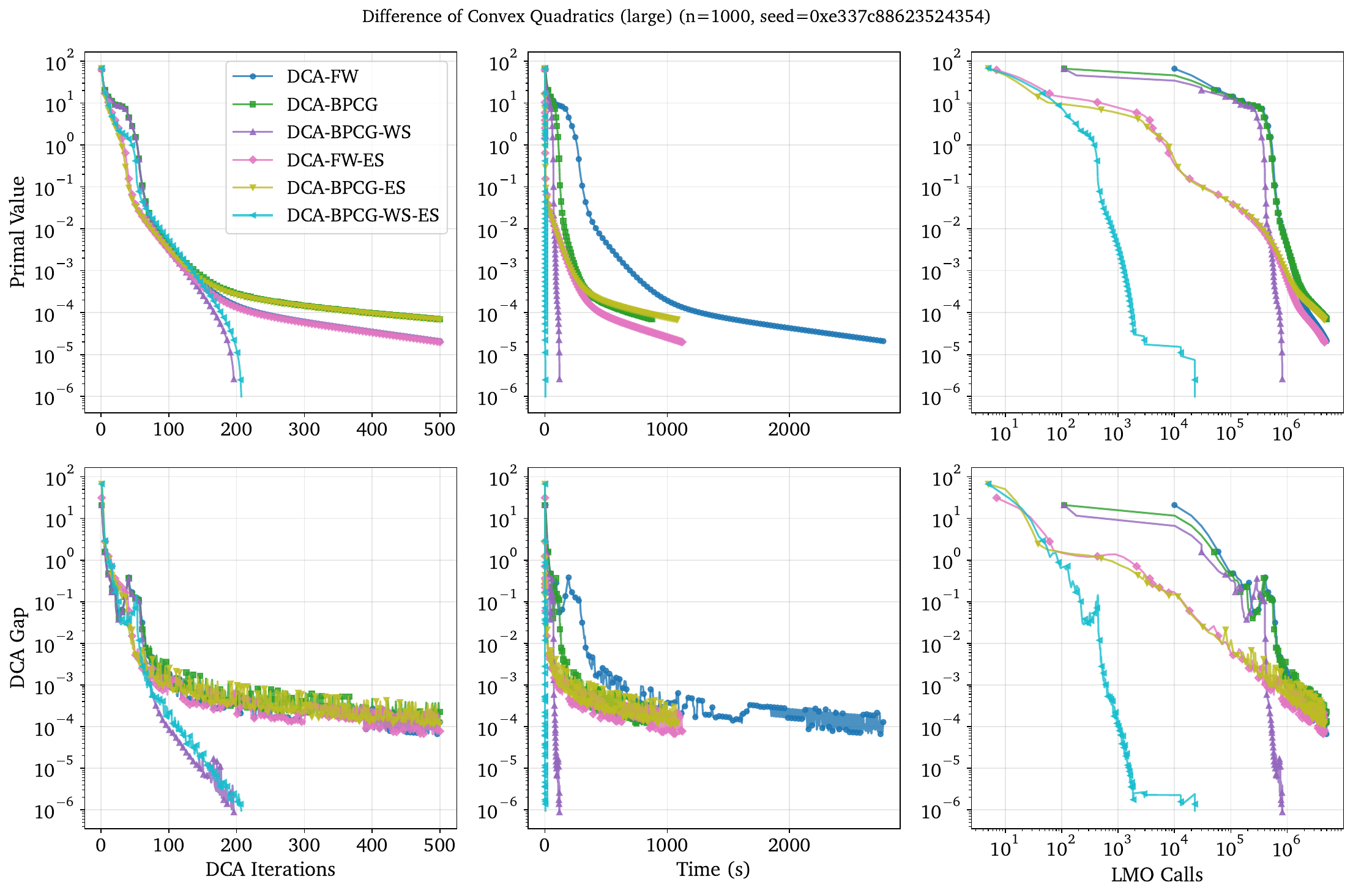}
    \caption{Run of the different variants on a Difference of Convex Quadratics example with seed $e337c88623524354$ with size $n = 1000$ from \cref{sec:diff-quadratics}.}
    \label{fig:diff-quadratics-1000-2}
\end{figure}

\begin{figure}[htbp]
    \centering
        \includegraphics[width=0.9\textwidth]{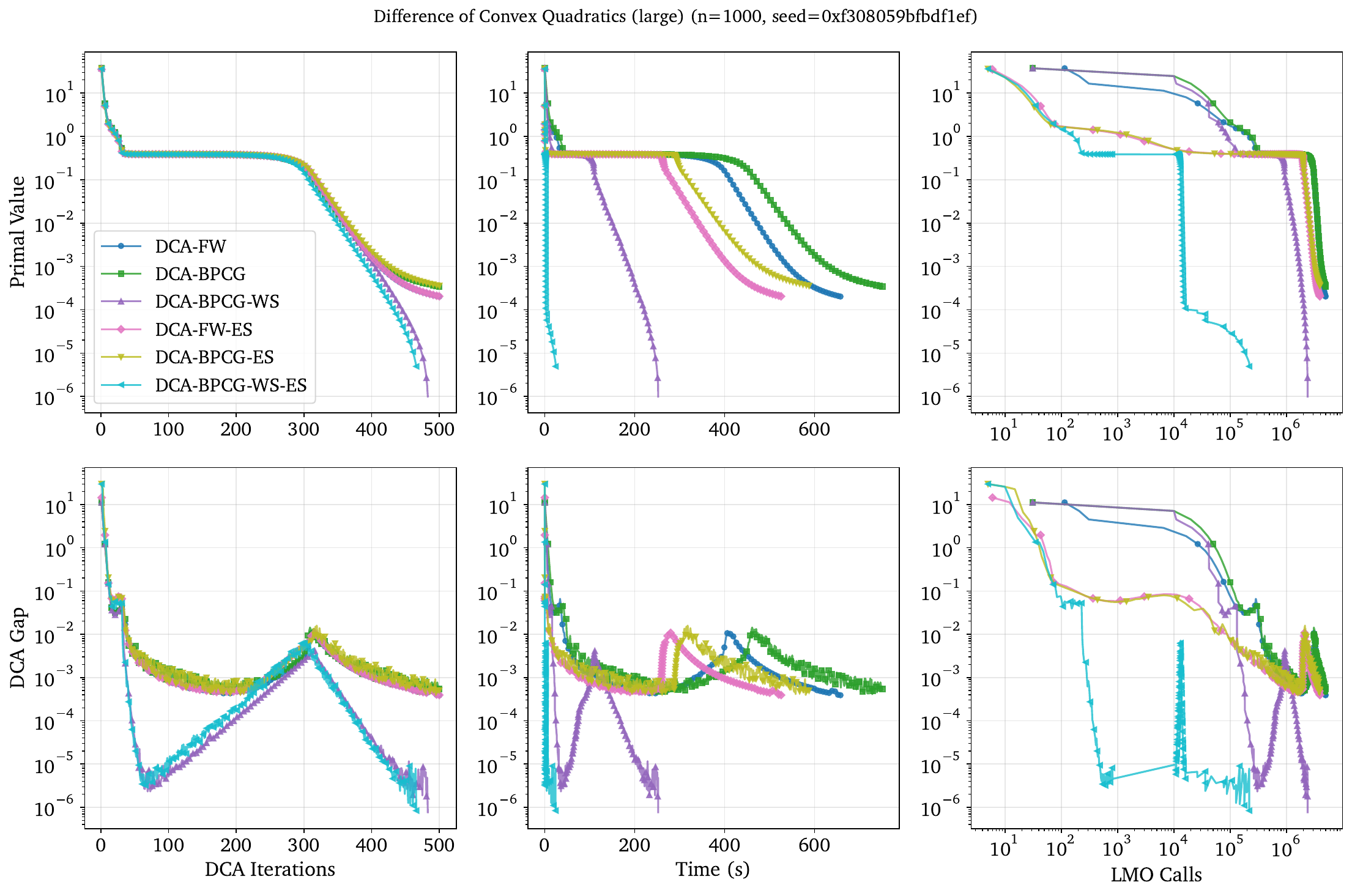}
    \caption{Run of the different variants on a Difference of Convex Quadratics example with seed $f308059bfbdf1ef$ with size $n = 1000$ from \cref{sec:diff-quadratics}.}
    \label{fig:diff-quadratics-1000-3}
\end{figure}

\begin{figure}[htbp]
    \centering
        \includegraphics[width=0.9\textwidth]{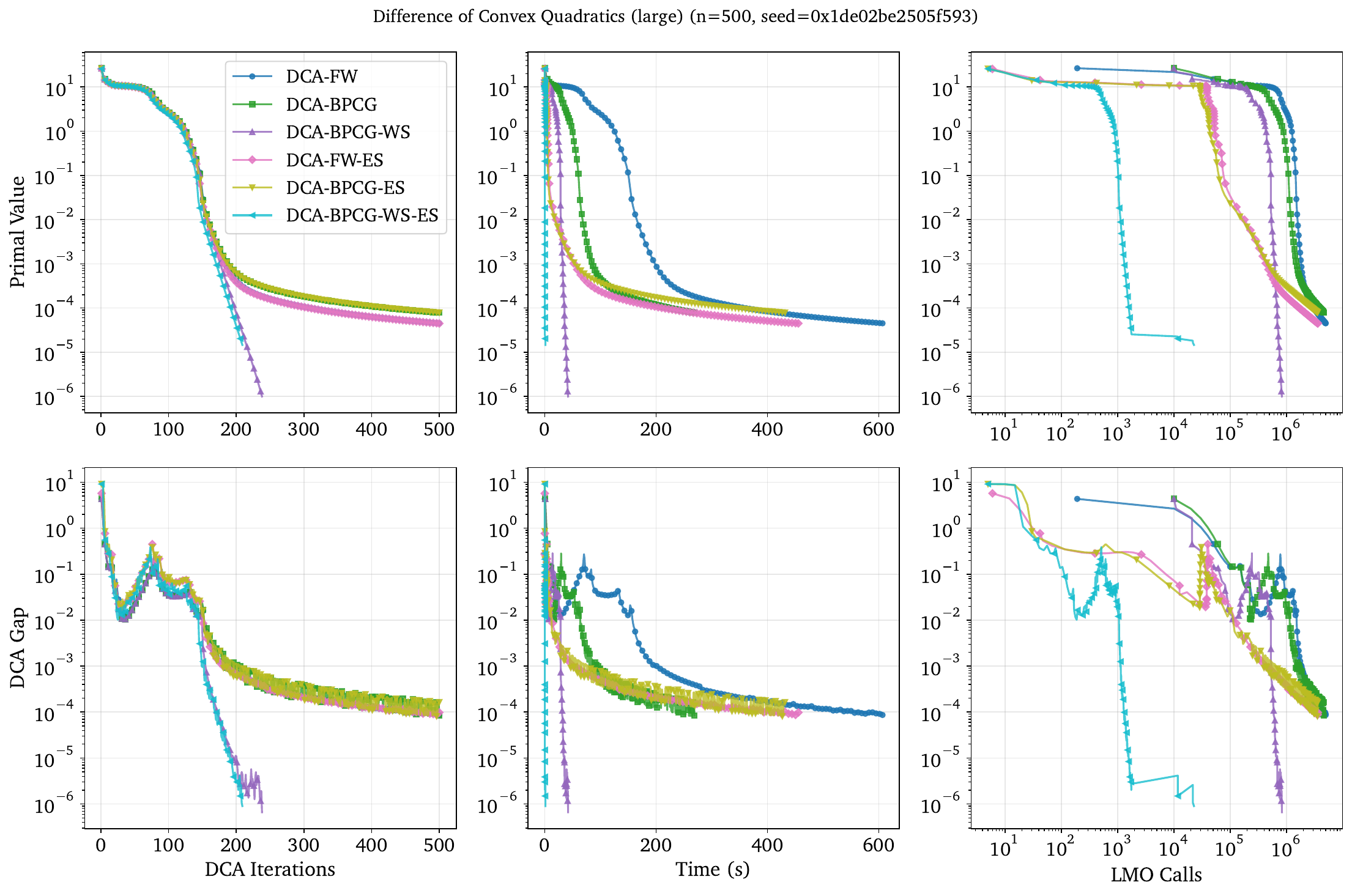}
    \caption{Run of the different variants on a Difference of Convex Quadratics example with seed $1de02be2505f593$ with size $n = 500$ from \cref{sec:diff-quadratics}.}
    \label{fig:diff-quadratics-500-1}
\end{figure}

In this section we present runs of some selected Difference of Convex Quadratic Functions examples from \cref{sec:diff-quadratics} in \cref{fig:diff-quadratics-1000-1,fig:diff-quadratics-1000-2,fig:diff-quadratics-1000-3,fig:diff-quadratics-500-1}; we took the largest sizes for presentation here as the results are more interesting.

\subsection{Hard DC function example}

\begin{figure}[htbp]
    \centering
        \includegraphics[width=0.9\textwidth]{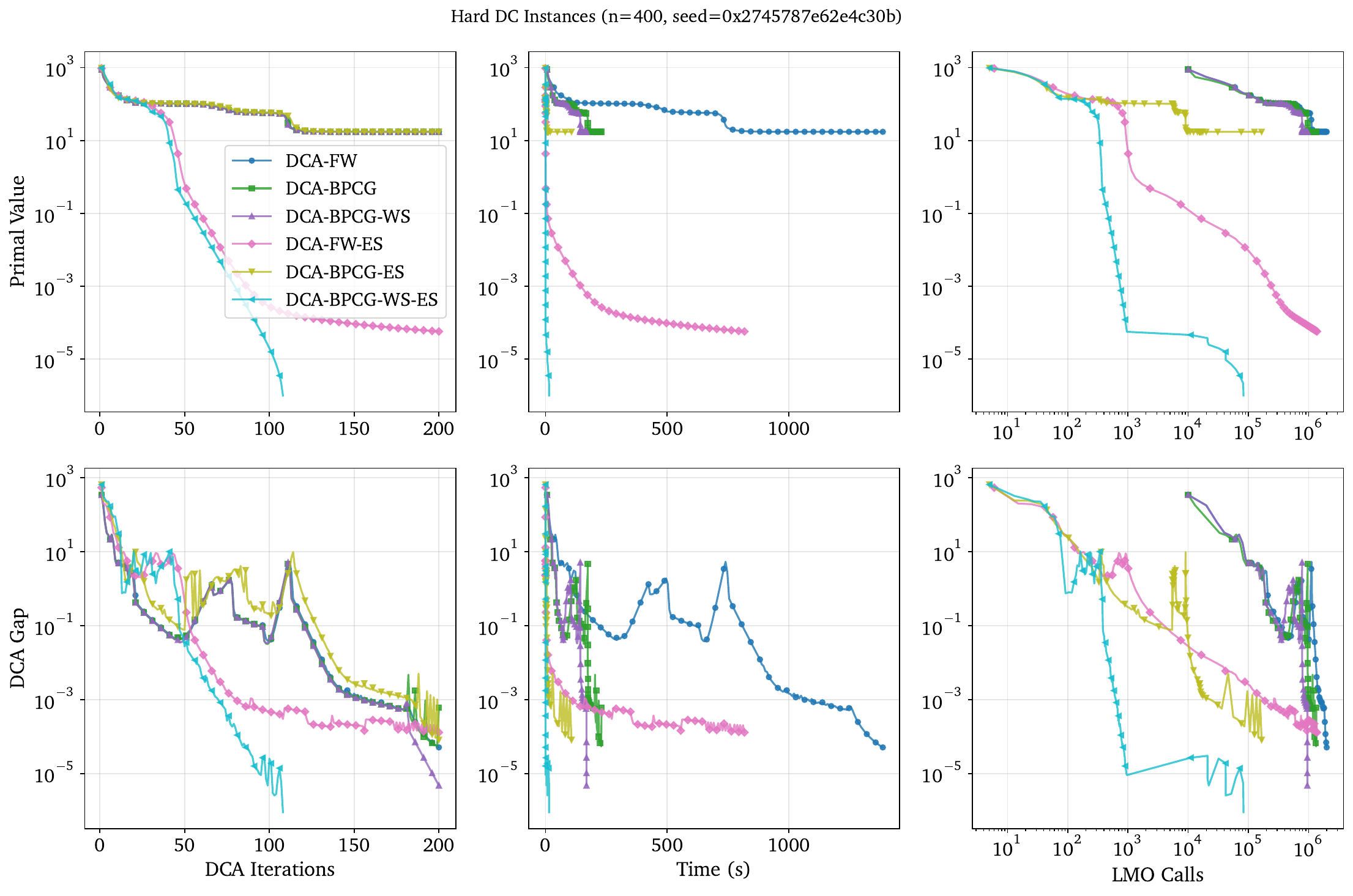}
    \caption{Run of the different variants on the hard DC function example with size $n = 400$ from \cref{sec:hard-dc}.}
    \label{fig:hard-dc}
\end{figure}

\begin{figure}[htbp]
    \centering
        \includegraphics[width=0.9\textwidth]{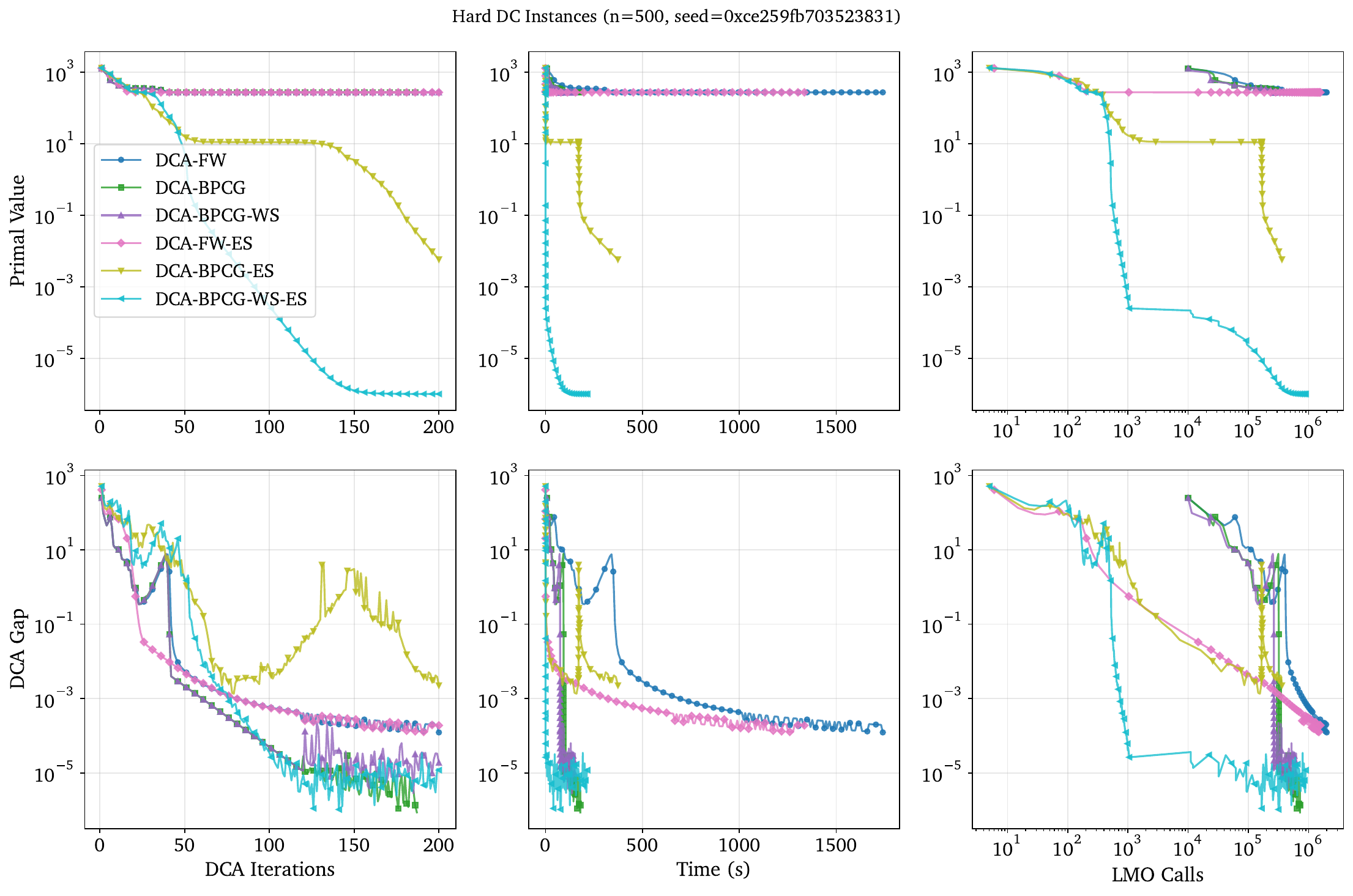}
    \caption{Run of the different variants on the hard DC function example with size $n = 500$ from \cref{sec:hard-dc}.}
    \label{fig:hard-dc-500-1}
\end{figure}

\begin{figure}[htbp]
    \centering
        \includegraphics[width=0.9\textwidth]{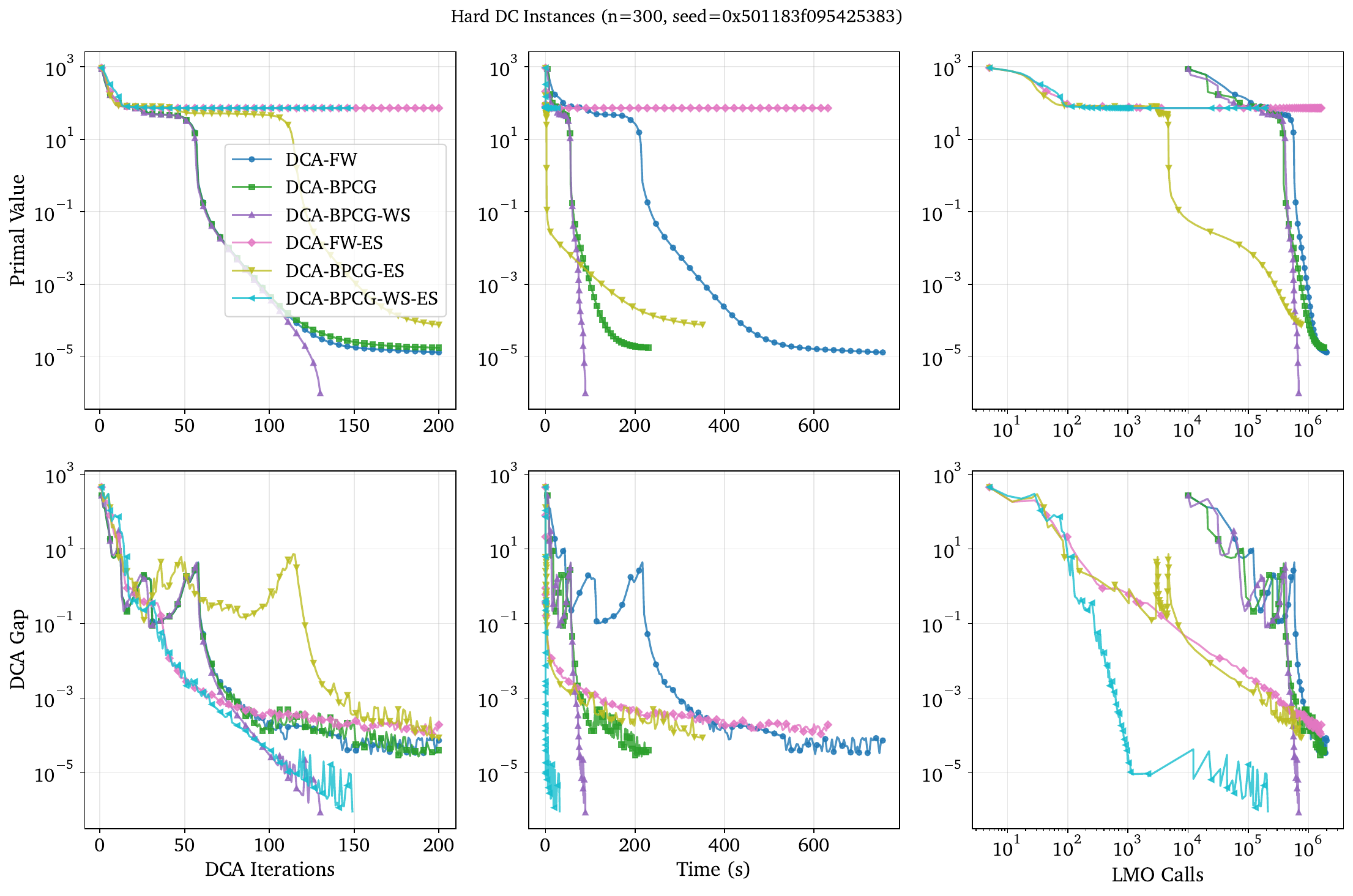}
    \caption{Run of the different variants on the hard DC function example with size $n = 300$ from \cref{sec:hard-dc}.}
    \label{fig:hard-dc-300-1}
\end{figure}

\begin{figure}[htbp]
    \centering
        \includegraphics[width=0.9\textwidth]{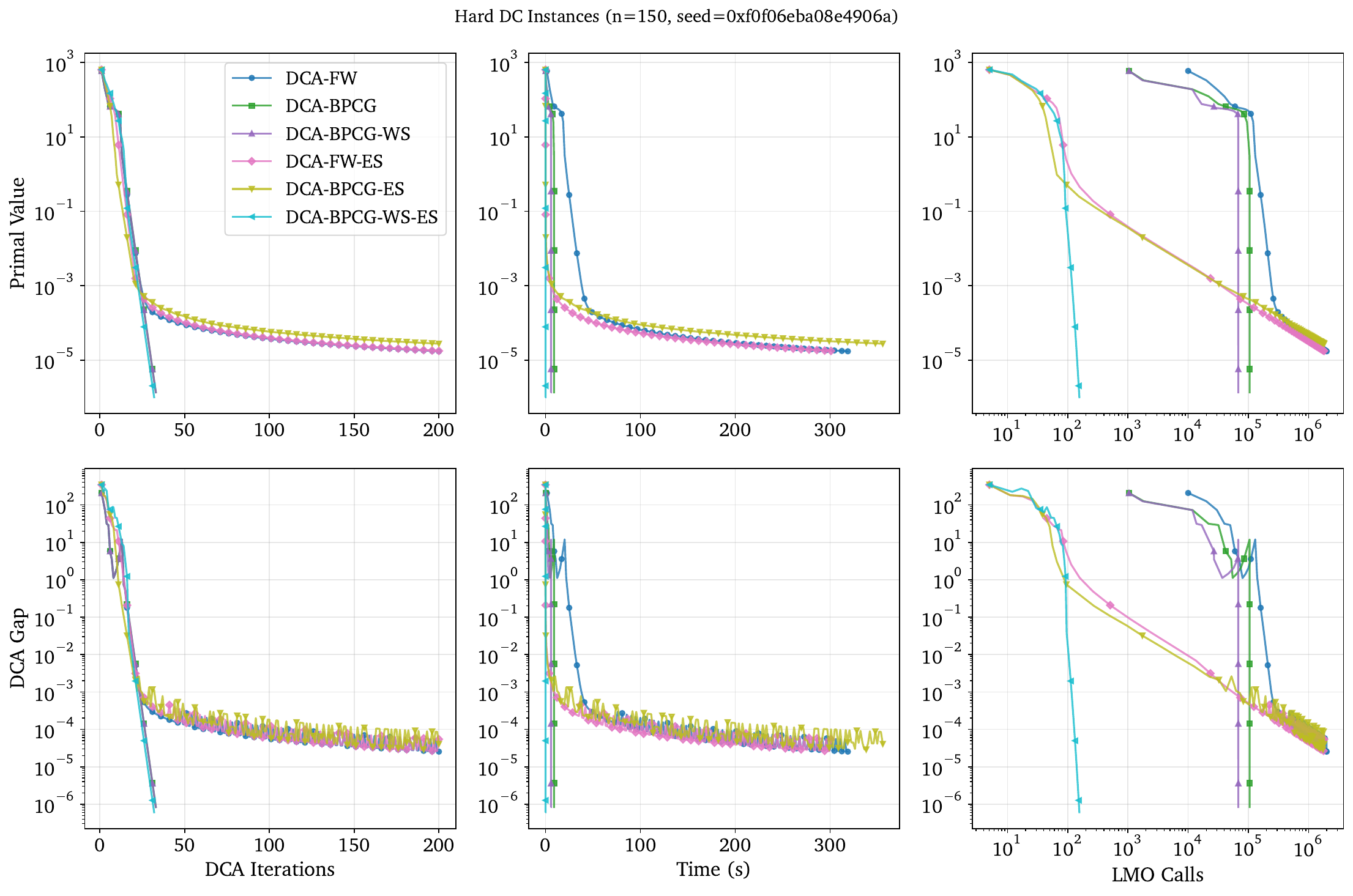}
    \caption{Run of the different variants on the hard DC function example with size $n = 150$ from \cref{sec:hard-dc}.}
    \label{fig:hard-dc-150-1}
\end{figure}

In this section we present runs of the hard DC function example from \cref{sec:hard-dc} in \cref{fig:hard-dc,fig:hard-dc-500-1,fig:hard-dc-300-1,fig:hard-dc-150-1}.

\subsection{Quadratic Assignment Problem}

\begin{figure}[htbp]
    \centering
        \includegraphics[width=0.9\textwidth]{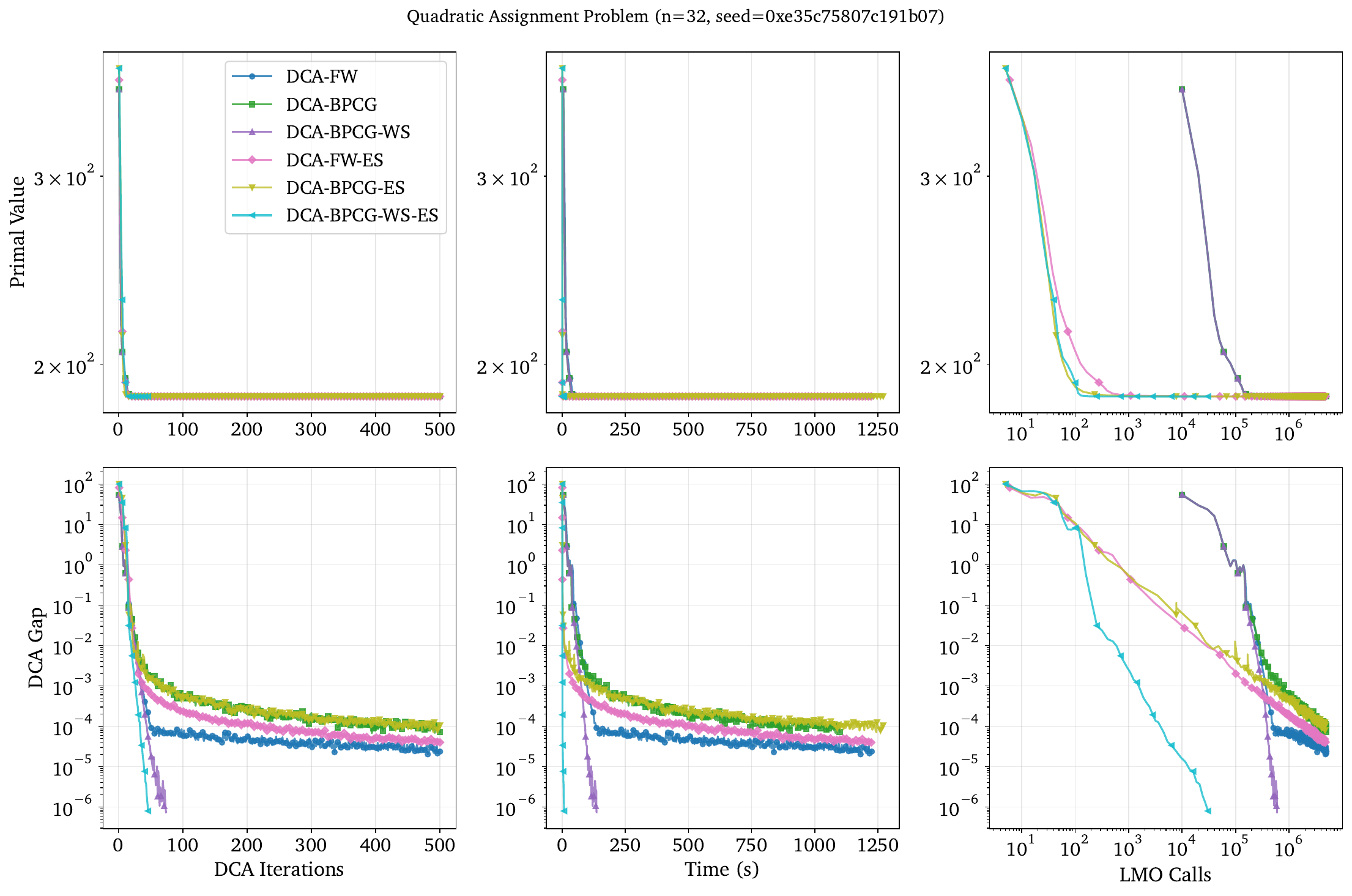}
    \caption{Run of the different variants on the small-sized\texttt{esc32b} QAP instance from \cref{sec:qap-instances}.}
    \label{fig:qap-esc32b}
\end{figure}

\begin{figure}[htbp]
    \centering
        \includegraphics[width=0.9\textwidth]{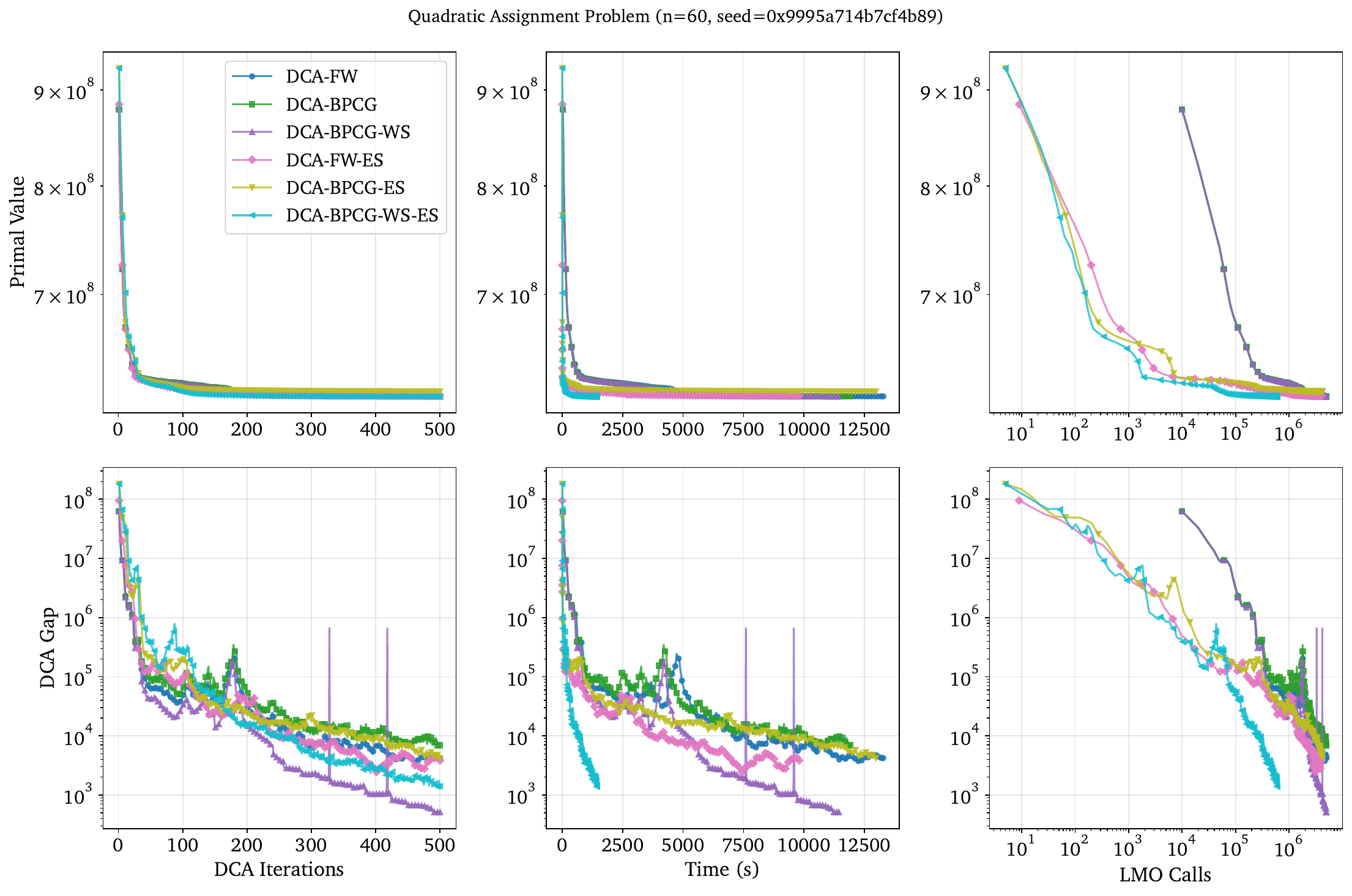}
    \caption{Run of the different variants on the medium-sized \texttt{tai60b} QAP instance from \cref{sec:qap-instances}.}
    \label{fig:qap-tai60b}
\end{figure}

\begin{figure}[htbp]
    \centering
        \includegraphics[width=0.9\textwidth]{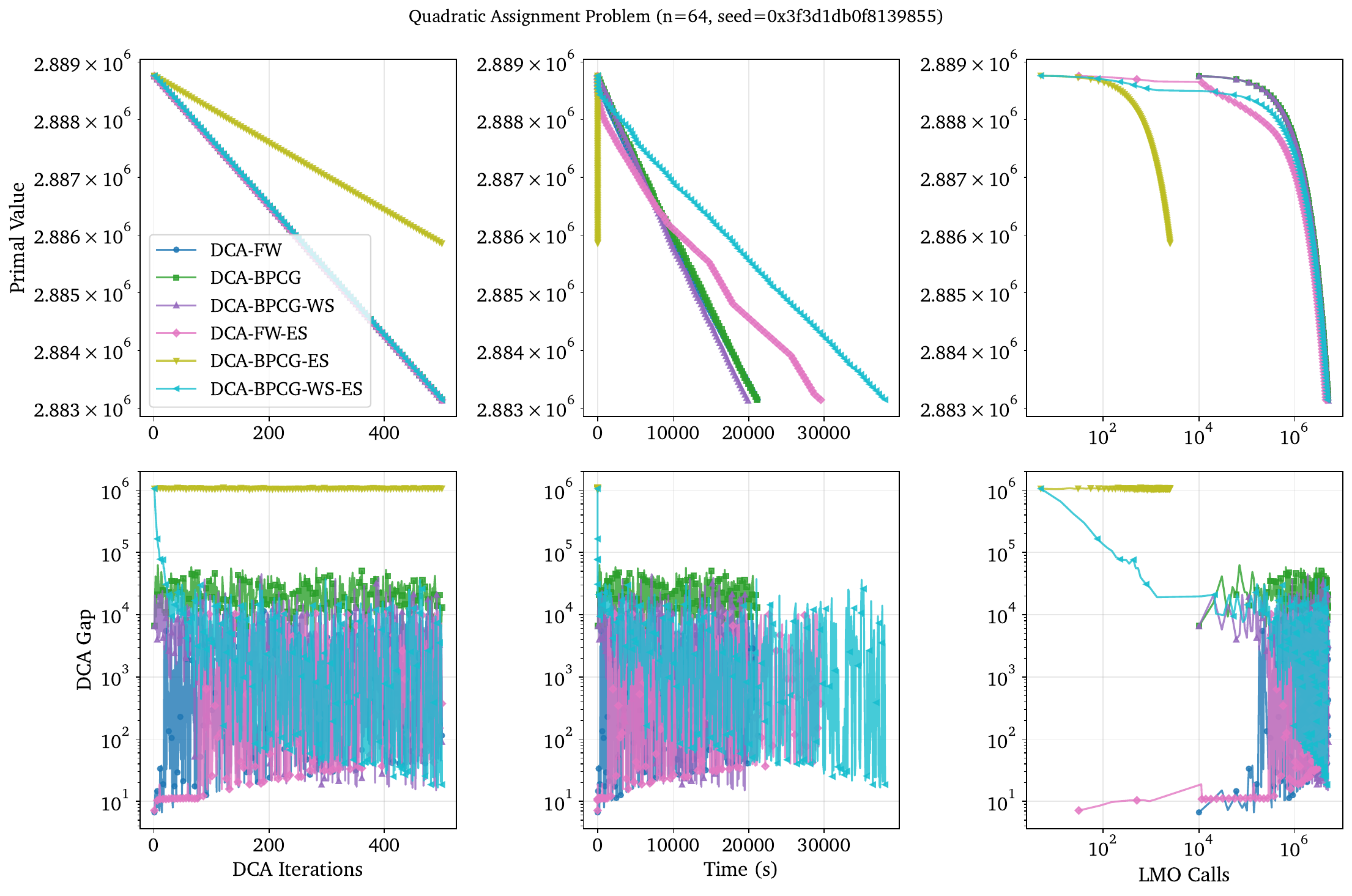}
    \caption{Run of the different variants on the medium-sized \texttt{tai64c} QAP instance from \cref{sec:qap-instances}. For this instance the DCA gaps are oscillating wildly.}
    \label{fig:qap-tai64c}
\end{figure}

\begin{figure}[htbp]
    \centering
        \includegraphics[width=0.9\textwidth]{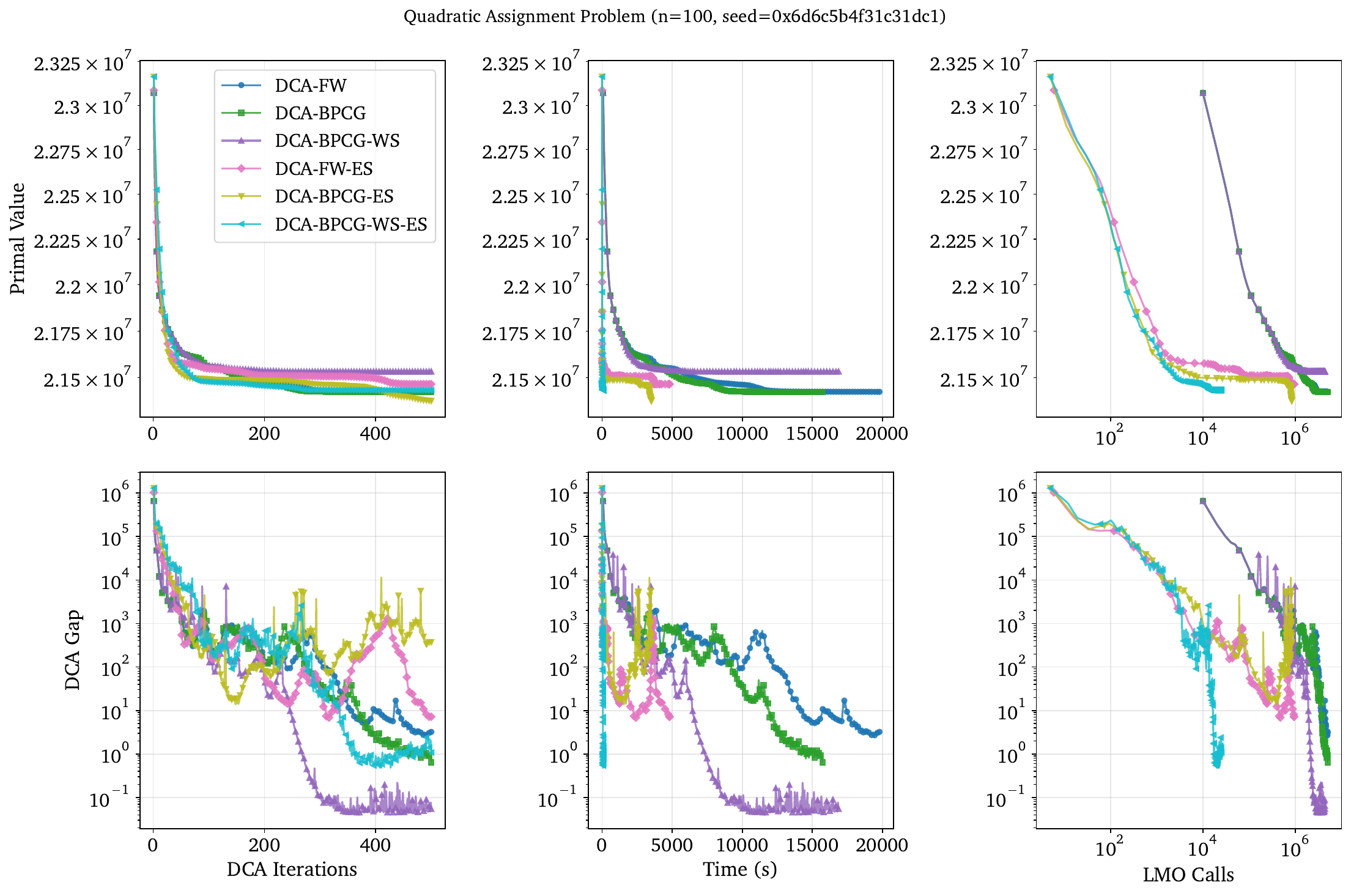}
    \caption{Run of the different variants on the large-sized \texttt{tai100a} QAP instance from \cref{sec:qap-instances}.}
    \label{fig:qap-tai100a}
\end{figure}

In this section we present runs of some selected QAP instances from \cref{sec:qap-instances} in \cref{fig:qap-esc32b,fig:qap-tai60b,fig:qap-tai64c,fig:qap-tai100a}.

\end{document}